\documentclass[a4paper,11pt]{scrartcl}
\usepackage{amsmath,amssymb}
\usepackage{amsthm}
\usepackage[square,sort&compress,numbers]{natbib}
\usepackage{graphicx}
\usepackage{url}
\usepackage{type1cm}
\usepackage{mathtools}\mathtoolsset{showonlyrefs=false}
\usepackage{algpseudocode,algorithm}
\usepackage{subfig}

\newcommand{\Authornote}{\renewcommand{\thefootnote}{\fnsymbol{footnote}}}
\newcommand{\authornote}{\Authornote\footnote}

\newcommand{\refalg}[1]{Algorithm~\ref{#1}}

\newcommand{\reffig}[1]{Figure~\ref{#1}}

\theoremstyle{remark}
\newtheorem{remark}{Remark}

\newcommand{\finbox}{\nolinebreak\hfill{\small $\blacksquare$}}

\newcommand{\MIN}{\mathop{\mathrm{Minimize}}}

\newcommand{\ST}{\mathop{\mathrm{subject~to}}}

\newcommand{\argmin}{\operatornamewithlimits{\mathrm{arg\,min}}}

\newcommand{\overRe}{\ensuremath{\Re\cup\{+\infty\}}}

\newcommand{\prox}{\mathop{\boldsymbol{\mathsf{prox}}}\nolimits}

\renewcommand{\Re}{\ensuremath{\mathbb{R}}}

\newcommand{\bi}[1]{\ensuremath{\boldsymbol{#1}}}

\newcommand{\rr}[1]{\ensuremath{\mathrm{#1}}}

\newcommand{\KC}{\ensuremath{\mathcal{K}}}
\newcommand{\LC}{\ensuremath{\mathcal{L}}}

\begin{document}

\begin{center}
  {\Large\bfseries\sffamily%
  Primal-dual algorithm for quasi-static contact problem }\\
  \medskip
  {\Large\bfseries\sffamily%
  with Coulomb's friction }%
  \par%
  \bigskip%
  {
  Yoshihiro Kanno~\authornote[2]{%
    Mathematics and Informatics Center, 
    The University of Tokyo, 
    Hongo 7-3-1, Tokyo 113-8656, Japan.
    E-mail: \texttt{kanno@mist.i.u-tokyo.ac.jp}. 
    }
  }
\end{center}

\begin{abstract}
  This paper presents a fast first-order method for solving the 
  quasi-static contact problem with the Coulomb friction. 
  It is known that this problem can be formulated as a second-order cone 
  linear complementarity problem, for which regularized or semi-smooth 
  Newton methods are widely used. 
  As an alternative approach, this paper develops a 
  method based on an accelerated primal-dual algorithm. 
  The proposed method is easy to implement, as most of computation 
  consists of additions and multiplications of vectors and matrices. 
  Numerical experiments demonstrate that this method outperforms a 
  regularized and smoothed Newton method for second-order cone 
  complementarity problems. 
\end{abstract}

\begin{quote}
  \textbf{Keywords}
  \par
  Optimization; 
  complementarity problem; 
  second-order cone; 
  primal-dual algorithm; 
  contact mechanics; 
  Coulomb friction. 
\end{quote}

\section{Introduction}

Frictional contact is ubiquitous in engineering applications, 
and its computational aspect has been studied 
extensively \citep{Wri06,AB08,Bro16,KO88}. 
Consider two bodies in the three-dimensional space. 
When the signed distance between their boundaries, 
called the {\em gap\/}, 
is positive, there exists no interaction in terms of forces. 
Alternatively, if the gap is equal to zero, then 
the contact pressure force, called the 
{\em normal reaction\/}, can be present at interface. 
This relation, which evidently possesses disjunction nature, is called the 
normal contact law. 
When the gap is equal to zero, 
the {\em tangential reaction\/} can also be present due to friction. 
The most fundamental model of friction is Coulomb's law. 
The relative tangential velocity of the two bodies at contact interface 
is whether equal to zero (said to be {\em stick\/}) 
or not (said to be {\em slip\/} or {\em slide\/}). 
The slip can occur when the magnitude of 
the tangential reaction is equal to a 
threshold value (which is proportional to the magnitude of the normal 
reaction), where the tangential reaction is in parallel with the 
relative tangential velocity in the opposite direction. 
In contrast, only stick is admissible when the magnitude of 
the tangential reaction is smaller than the threshold. 
Thus, the tangential contact law also possesses disjunction nature. 

The disjunction nature in frictional contact can be treated within the 
framework of complementarity problem. 
To date, diverse formulations, as well as diverse numerical methods, 
have been proposed; see, e.g., \citet{AB08} and \citet{Wri06}. 
This reflects the fact that, as \citet{ABH18} concluded through 
comprehensive numerical experiments, ``there is no universal solver'' 
for the frictional contact problem. 
Indeed, to the best of the author's knowledge, there is no algorithm 
that has guarantee of convergence for the three-dimensional quasi-static 
incremental contact problem with the Coulomb friction. 

In this paper, we address the three-dimensional quasi-static incremental 
problem of elastic bodies subjected to the unilateral contact with 
the Coulomb friction, under the assumption of small deformation and 
linear elasticity. 
For solving this problem, this paper presents an algorithm based on 
the {\em primal-dual algorithm\/} 
(a.k.a.\ the primal-dual hybrid gradient algorithm) 
\citep{CP11,CP16b,HYY14,MP18}. 
The primal-dual algorithm has received considerable attention in 
applications to large-scale optimization problems arising in image 
processing \citep{CP11,CP16,BR12,EZC10,CLO14}. 

If we approximate the friction cone (see \eqref{eq:friction.cone} in 
section~\ref{sec:problem.friction} for definition of the friction cone) 
as a polyhedral cone, 
the problem considered in this paper is reduced to 
a {\em linear complementarity problem\/} \citep{Kla86}. 
In contrast, without use of this approximation, the problem is 
formulated as a {\em nonlinear complementarity problem\/}. 
Regularized or semi-smooth Newton methods are mainly used to solve such 
formulations \citep{AC91,Chr02,CKPS98,Ren13,APRQDB14}. 
\citet{KMPc06} showed that the problem can be recast as 
a {\em second-order cone linear complementarity problem\/}. 
To this formulation, various methods developed in the field of 
mathematical optimization may be applicable (although existing methods 
do not have a proof of convergence for this formulation, 
due to the lack of monotonicity; see the formulation in 
\citet[section~4]{KMPc06}). 
See \citet{Yos12} and \citet{CP12} for survey on the second-order cone 
complementarity problem and its numerical solutions. 
For example, a regularized smoothing Newton method proposed by 
\citet{HYF05} was adopted in \citet{KMPc06}. 

Recently, {\em accelerated gradient methods\/} 
(a.k.a.\ {\em optimal first-order methods\/}) have been successfully 
developed for solving problems, especially large-scale ones, in applied 
and computational mechanics. 
Namely, {\em accelerated proximal gradient methods\/} have been proposed 
for solving the incremental problems in elastoplasticity with various 
yield criteria \citep{Kan16,Kan20,SK18,SK20}, as well as the bi-modulus 
elasticity problem \citep{Kan21}. 
It is worth noting that the problems dealt with in the literature 
above are formulated as convex optimization problems. 
Specifically, the problem in \citep{Kan16} can be recast 
as a {\em quadratic programming\/} (QP) problem, 
the one in \citep{SK18} can be recast as 
a {\em second-order cone programming\/} (SOCP) problem, 
and the ones in \citep{SK20,Kan21} can be recast as 
{\em semidefinite programming\/} (SDP) problems. 
The numerical experiments demonstrate that these accelerated proximal 
gradient methods outperform standard solvers implementing primal-dual 
interior-point methods for QP, SOCP, and SDP. 
As other types of accelerated gradient methods in mechanics, the reader 
may refer to an accelerated Uzawa method for the 
frictionless contact problem \citep{Kan20Uzawa} and 
an accelerated steepest descent method for the elasticity problem of 
trusses with material nonlinearity \citep{FK19}. 

In the field of computer graphics, \citet{MHNT15} proposed an 
{\em accelerated projected gradient method\/} to (approximately) solve a 
problem stemming from time-discretization of a rigid multi-body dynamical 
system involving frictional contact. 
Through numerical experiments, \citet{MFJN17} concluded that this method 
is most efficient compared with conventional first-order methods (the 
projected Jacobi method and the projected Gauss--Seidel method) in this 
research area and second-order optimization methods 
(primal-dual interior-point methods). 
It should be clear that the method of \citet{MHNT15} does {\em not\/} 
solve the problem to be solved in dynamic simulation, but is a method to 
solve a modified problem which is easier than the original one. 
More concretely, the original problem is a nonlinear complementarity 
problem that does {\em not\/} correspond to the optimality condition of 
an optimization problem, while \citet{MHNT15} solves a convex 
optimization problem that is obtained by adding modification to the 
original problem. 
Such artificial modification certainly alters physical phenomena, and 
causes artifacts in simulation results as illustrated in 
\citet[section~2.3]{MHNT15}. 
It is worth noting that there exists some numerical methods for solving 
the original problem; e.g., 
nonsmooth Newton methods \citep{Cad09,BCDA11}, 
a fixed-point method that sequentially solves convex optimization 
problems \citep{ACLM11}, etc. 
These methods are not first-order methods. 
In contrast, in this paper we consider the quasi-static incremental 
problem of elastic bodies, and attempt to solve the problem itself, i.e., 
without adding any modification, via an accelerated first-order method.


The paper is organized as follows. 
Section~\ref{sec:problem} summarizes the frictional contact problem 
considered in this paper. 
As a main contribution, section~\ref{sec:algorithm} presents an 
algorithm for solving this problem, based on 
an accelerated primal-dual algorithm. 
Section~\ref{sec:ex} performs numerical experiments. 
Section~\ref{sec:conclusion} presents some conclusions.


In our notation, 
${}^{\top}$ denotes the transpose of a vector or a matrix. 
For $\bi{x} \in \Re^{n}$, we use $\| \bi{x} \|$ to denote its 
Euclidean norm, i.e., $\| \bi{x} \| = \sqrt{\bi{x}^{\top} \bi{x}}$. 
For $\bi{x}$, $\bi{y} \in \Re^{n}$, 
we write $\bi{x} \perp \bi{y}$ if $\bi{x}^{\top} \bi{y} = 0$. 
We use $\langle \bi{x},\bi{y} \rangle$ to denote 
$\bi{x}^{\top} \bi{y}$. 
For a closed convex function $f : \Re^{n} \to \overRe$, 
its {\em proximal mapping\/} is defined by 
\begin{align}
  \prox_{f}(\bi{x}) 
  = \argmin_{\bi{z} \in \Re^{n}} \Bigl\{
  f(\bi{z}) + \frac{1}{2} \| \bi{z} - \bi{x} \|^{2}
  \Bigr\} . 
  \label{eq:def.proximal.mapping}
\end{align}
For $S \subseteq \Re^{n}$, we use 
$\delta_{S} : \Re^{n} \to \overRe$ to denote its 
{\em indicator function\/}, i.e., 
\begin{align*}
  \delta_{S}(\bi{x}) = 
  \begin{dcases*}
    0 
    & if $\bi{x} \in S$, \\
    +\infty 
    & otherwise. 
  \end{dcases*}
\end{align*}
We use $\Pi_{S}(\bi{x}) \in \Re^{n}$ to denote the {\em projection\/} 
of $\bi{x} \in \Re^{n}$ onto $S$, i.e., 
\begin{align*}
  \Pi_{S}(\bi{x}) 
  = \argmin_{\bi{z} \in S} 
  \{ \| \bi{z} - \bi{x} \| \} . 
\end{align*}
For a nonempty convex cone $C \subseteq \Re^{n}$, 
define the {\em dual cone\/} by 
\begin{align*}
  C^{*} 
  = \{ \bi{s} \in \Re^{n} 
  \mid
  \langle \bi{s} , \bi{x} \rangle \ge 0 
  \
  (\forall \bi{x} \in C)
  \} . 
\end{align*}
We readily see that 
\begin{align}
  \inf_{\bi{x} \in \Re^{n}} \{ 
  \langle \bi{s} , \bi{x} \rangle + \delta_{C}(\bi{x}) \}
  = -\delta_{C^{*}}(\bi{s}) 
  \label{eq:dual.cone.1}
\end{align}
holds.

\section{Fundamentals: Quasi-static incremental problem}
\label{sec:problem}

This section summarizes the quasi-static incremental analysis of an 
elastic solid, with the Coulomb friction for the unilateral contact. 
For fundamentals of contact mechanics, see, e.g., 
\citet{DL76} and \citet{Wri06}. 

\subsection{Contact kinematics}

Consider an elastic body in the three-dimensional space. 
The body is discretized according to the conventional finite element 
method. 
Let $d$ denote the number of degrees of freedom of the nodal 
displacements. 
We use $\bi{u} \in \Re^{d}$ to denote the nodal displacement vector. 

To investigate the time evolution in the specified time interval $[0,T]$, 
suppose that the time interval is subdivided into finitely many 
intervals. 
In the {\em quasi-static analysis\/}, we assume that the inertia term in 
the equation of motion is negligible. 
This assumption is applicable if the external force applied to the 
elastic body changes sufficiently slowly. 
For a specific time subinterval, denoted by $[t^{l},t^{l+1}]$, let 
$\bi{u}^{l}$ and $\bi{u}^{l+1}$ denote the nodal displacements at 
time $t^{l}$ and $t^{l+1}$, respectively. 
The {\em incremental problem\/} solved for time $t^{l+1}$ is to find 
$\bi{u}^{l+1}$, or, equivalently, to find 
the {\em incremental displacement\/} \citep{Kla99}
\begin{align*}
  \Delta\bi{u} 
  = \bi{u}^{l+1} - \bi{u}^{l} . 
\end{align*}

\begin{figure*}[tbp]
  \centering
  \includegraphics[scale=0.60]{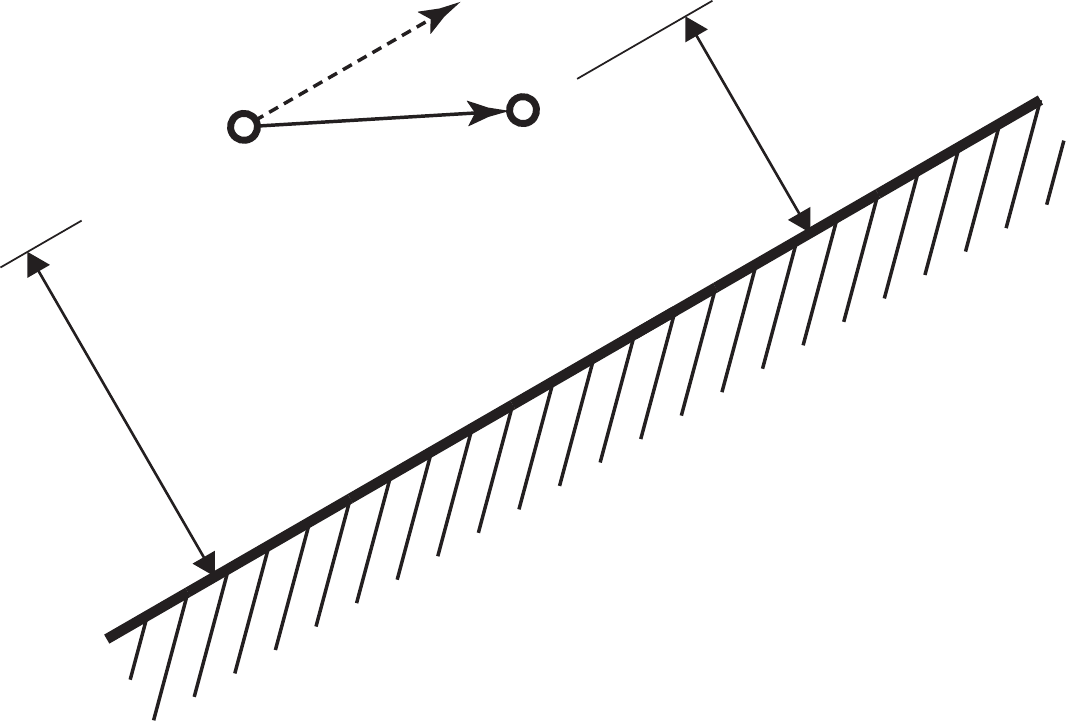}
  \begin{picture}(0,0)
    \put(-200,-100){
    \put(49,218){{\small $\Delta\bi{u}_{\rr{t}j}$}}
    \put(142,207){{\small $\hat{g}_{j}(\Delta\bi{u})$}}
    \put(21,150){{\small $g_{j}$}}
    \put(25,200){{\small $t=t^{l}$}}
    \put(90,195){{\small $t=t^{l+1}$}}
    \put(96,128){{\small obstacle}}
    }
  \end{picture}
  \caption[]{Contact candidate node and obstacle. }
  \label{fig:local_coordinate}
\end{figure*}

For simplicity, we restrict ourselves to the case that the boundary of an 
elastic body can possibly touch the surface of a fixed rigid obstacle; 
extension of the proposed method to the case that two elastic bodies 
can touch each other is straightforward. 
We assume that a set of {\em contact candidate nodes\/} (i.e., nodes 
that can possibly contact with the obstacle at time $t=t^{l+1}$) is 
specified a priori. 
\reffig{fig:local_coordinate} depicts one of contact candidate nodes. 
It should be clear that we do not know in advance whether each contact 
candidate node contacts with the obstacle or not at $t=t^{l+1}$ 
(because $\Delta\bi{u}$ is unknown). 

Let $c$ denote the number of contact candidate nodes. 
For contact candidate node $j$ $(j=1,\dots,c)$, 
let $g_{j}$ denote the {\em initial gap\/} (i.e., 
the distance between node $j$ and the obstacle at time $t=t^{l}$), 
which is a nonnegative constant; see \reffig{fig:local_coordinate}. 
The gap at time $t=t^{l+1}$, denoted by $\hat{g}_{j}(\Delta\bi{u})$, is 
given in the form 
\begin{align}
  \hat{g}_{j}(\Delta\bi{u}) 
  = g_{j} - \bi{t}_{\rr{n}j}^{\top} \, \Delta\bi{u} 
  \label{eq:kinematics.1}
\end{align}
with a constant vector $\bi{t}_{\rr{n}j} \in \Re^{d}$, 
because we assume small deformation. 
The incremental displacement of node $j$ can be decomposed additively 
into two components, which are normal and tangential to the obstacle 
surface. 
The tangential component, denoted by $\Delta\bi{u}_{\rr{t}j}$, is given 
in the form 
\begin{align}
  \Delta\bi{u}_{\rr{t}j} 
  = T_{\rr{t}j}^{\top} \, \Delta\bi{u} , 
  \label{eq:kinematics.2}
\end{align}
where $T_{\rr{t}j} \in \Re^{d \times 2}$ is a constant matrix; 
see \citet{Kla99} for more details. 
For notational simplicity, define $T_{\rr{n}} \in \Re^{d \times c}$ 
and $T_{\rr{t}} \in \Re^{d \times 2c}$ by 
\begin{align*}
  T_{\rr{n}} = 
  \begin{bmatrix}
    \bi{t}_{\rr{n}1} & \cdots & \bi{t}_{\rr{n}c} \\
  \end{bmatrix}
  , \quad
  T_{\rr{t}} = 
  \begin{bmatrix}
    T_{\rr{t}1} & \cdots & T_{\rr{t}c} \\
  \end{bmatrix}
  . 
\end{align*}

If the gap $\hat{g}_{j}(\Delta\bi{u})$ is equal to zero, then 
the contact reaction force can be present at node $j$. 
The reaction also has two components, denoted by 
$r_{\rr{n}j} \in \Re$ and $\bi{r}_{\rr{t}j} \in \Re^{2}$, which are 
normal and tangential to the obstacle surface, respectively. 
It is worth noting that $r_{\rr{n}j}$ stems from a kinematic constraint 
that node $j$ cannot penetrate the obstacle, 
while $\bi{r}_{\rr{t}j}$ is due to friction. 
For notational convenience, define $\bi{r}_{j} \in \Re^{3}$, 
$\bi{r} \in \Re^{3c}$, $\bi{r}_{\rr{n}} \in \Re^{c}$, and 
$\bi{r}_{\rr{t}} \in \Re^{2c}$ by 
\begin{align*}
  \bi{r}_{j} = 
  \begin{bmatrix}
    r_{\rr{n}j}   \\  \bi{r}_{\rr{t}j} \\
  \end{bmatrix}
  , \quad
  \bi{r} = 
  \begin{bmatrix}
    \bi{r}_{1} \\ \vdots \\ \bi{r}_{c} \\
  \end{bmatrix}
  , \quad
  \bi{r}_{\rr{n}} = 
  \begin{bmatrix}
    r_{\rr{n}1} \\ \vdots \\ r_{\rr{n}c} \\
  \end{bmatrix}
  , \quad
  \bi{r}_{\rr{t}} = 
  \begin{bmatrix}
    \bi{r}_{\rr{t}1} \\ \vdots \\ \bi{r}_{\rr{t}c} \\
  \end{bmatrix}
  . 
\end{align*}

\subsection{Coulomb's friction law under unilateral contact}
\label{sec:problem.friction}

Let $\mu$ denote the coefficient of friction, which is a positive 
constant. 
For the incremental problem, the {\em Coulomb friction law\/} for 
unilateral contact is given by \citep{DL76,Wri06} 
\begin{subequations}\label{eq:Coulomb.law}%
  \begin{alignat}{3}
    & {-}\mu r_{\rr{n}j} \ge \| \bi{r}_{\rr{t}j} \| , 
    \label{eq:Coulomb.law.1} \\
    & {-}\mu r_{\rr{n}j} > \| \bi{r}_{\rr{t}j} \| 
    &&\quad\Rightarrow\quad
    \Delta\bi{u}_{\rr{t}j} = \bi{0} , 
    \label{eq:Coulomb.law.2} \\
    & {-}\mu r_{\rr{n}j} = \| \bi{r}_{\rr{t}j} \| > 0 
    &&\quad\Rightarrow\quad
    \exists \alpha \ge 0 : \
    \Delta\bi{u}_{\rr{t}j} = -\alpha \bi{r}_{\rr{t}j} . 
    \label{eq:Coulomb.law.3}
  \end{alignat}
\end{subequations}
Here, \eqref{eq:Coulomb.law.1} is inclusion of the reaction into the 
friction cone, which implies $r_{\rr{n}j} \le 0$, i.e., there exists 
{\em no adhesion\/}. 
Disjunction of the {\em sticking\/} and {\em slipping\/} states is 
described in \eqref{eq:Coulomb.law.2} and \eqref{eq:Coulomb.law.3}. 

Besides the non-adhesion condition (i.e., $r_{\rr{n}j} \le 0$), the 
{\em unilateral contact\/} condition consists of \citep{DL76,Wri06} 
\begin{subequations}\label{eq:unilateral.law}%
  \begin{alignat}{3}
    & \hat{g}_{j}(\Delta\bi{u}) \ge 0 , 
    \label{eq:unilateral.law.1} \\
    & \hat{g}_{j}(\Delta\bi{u}) > 0 
    &&\quad\Rightarrow\quad
    r_{\rr{n}j} = 0 , 
    \label{eq:unilateral.law.2} \\
    & r_{\rr{n}j} < 0 
    &&\quad\Rightarrow\quad
    \hat{g}_{j}(\Delta\bi{u}) = 0 . 
    \label{eq:unilateral.law.3}
  \end{alignat}
\end{subequations}
Here, \eqref{eq:unilateral.law.1} is the {\em non-penetration\/} condition. 
We say that node $j$ is in the {\em free\/} state and 
{\em in contact\/} (with nonzero reaction), respectively, if 
\eqref{eq:unilateral.law.2} and \eqref{eq:unilateral.law.3} hold. 

It is known that \eqref{eq:Coulomb.law} and \eqref{eq:unilateral.law} 
can be rewritten equivalently as \citep{MPcS02,MPc00} 
\begin{subequations}\label{eq:friction.law}%
  \begin{align}
    & \hat{g}_{j}(\Delta\bi{u}) \ge 0 ,  
    \label{eq:friction.law.1} \\
    & {-}\mu r_{\rr{n}j} \ge \| \bi{r}_{\rr{t}j} \| , \\
    & \langle \bi{r}_{\rr{t}j}, \Delta\bi{u}_{\rr{t}j} \rangle
    - \langle r_{\rr{n}j} , 
    \hat{g}_{j}(\bi{u}) + \mu \| \Delta\bi{u}_{\rr{t}j}\| \rangle = 0 . 
  \end{align}
\end{subequations}
Let $F \subset \Re^{3}$ denote the {\em friction cone\/}, i.e., 
\begin{align}
  F = \{
  (r_{\rr{n}},\bi{r}_{\rr{t}}) \in \Re \times \Re^{2} 
  \mid
  -\mu r_{\rr{n}} \ge \| \bi{r}_{\rr{t}} \|
  \} . 
  \label{eq:friction.cone}
\end{align}
The dual cone of $F$ is 
\begin{align*}
  F^{*} = \{
  (v_{\rr{n}},\bi{v}_{\rr{t}}) \in \Re \times \Re^{2}  
  \mid
  -v_{\rr{n}} \ge \mu \| \bi{v}_{\rr{t}} \|
  \} . 
\end{align*}
Since \eqref{eq:friction.law.1} holds if and only if 
$-\hat{g}_{j}(\Delta\bi{u}) - \mu \| \Delta\bi{u}_{\rr{t}j} \|  \in F^{*}$, 
we see that \eqref{eq:friction.law} can be recast as the following cone 
complementarity condition \citep[section~10.3.4.1]{Kan11}: 
\begin{align}
  F^{*} \ni 
  \begin{bmatrix}
    -\hat{g}_{j}(\Delta\bi{u}) - \mu \| \Delta\bi{u}_{\rr{t}j} \| \\
    \Delta\bi{u}_{\rr{t}j} \\
  \end{bmatrix}
  \perp
  \begin{bmatrix}
    r_{\rr{n}j} \\ \bi{r}_{\rr{t}j} \\
  \end{bmatrix}
  \in F . 
  \label{eq:complementarity.form.1}
\end{align}

\subsection{Equilibrium equation}

Let $K \in \Re^{d \times d}$ denote the stiffness matrix of the elastic 
body, which is symmetric and positive definite. 
We use $\bi{p}^{l} \in \Re^{d}$ to denote the nodal external force 
applied to the elastic body at time $t^{l}$. 

At an equilibrium state, 
the internal force, the external force, and the reaction are balanced. 
At time $t^{l}$, this balance law, called the 
{\em equilibrium equation\/}, is given by \citep{Kla99,Wri06} 
\begin{align}
  K \bi{u}^{l} - \bi{p}^{l} 
  = T_{\rr{n}} \bi{r}_{\rr{n}}^{l} + T_{\rr{t}} \bi{r}_{\rr{t}}^{l} . 
  \label{eq:equilibrium.eq.1}
\end{align}
In the incremental problem for time $t^{l+1}$, we have already known 
$\bi{u}^{l}$, $\bi{r}_{\rr{n}}^{l}$, and $\bi{r}_{\rr{t}}^{l}$ 
satisfying \eqref{eq:equilibrium.eq.1}, and attempt to find 
$\bi{u}^{l+1}$, $\bi{r}_{\rr{n}}^{l+1}$, and $\bi{r}_{\rr{t}}^{l+1}$ 
satisfying the equilibrium equation at time $t^{l+1}$, i.e., 
\begin{align}
  K (\bi{u}^{l} + \Delta\bi{u}) - \bi{p}^{l+1} 
  = T_{\rr{n}} \bi{r}_{\rr{n}}^{l+1} 
  + T_{\rr{t}} \bi{r}_{\rr{t}}^{l+1} . 
  \label{eq:equilibrium.eq.2}
\end{align}

For notational simplicity, we use $\bi{p} \in \Re^{d}$, 
$\bi{r}_{\rr{n}} \in \Re^{c}$, and $\bi{r}_{\rr{t}} \in \Re^{2c}$ to 
denote 
\begin{align*}
  \bi{p} = \bi{p}^{l+1} - K \bi{u}^{l} , 
  \quad
  \bi{r}_{\rr{n}} = \bi{r}_{\rr{n}}^{l+1} 
  \quad
  \bi{r}_{\rr{t}} = \bi{r}_{\rr{t}}^{l+1} . 
\end{align*}
From \eqref{eq:equilibrium.eq.1} and \eqref{eq:equilibrium.eq.2}, 
we have 
\begin{align}
  K \, \Delta\bi{u} - \bi{p} 
  = T_{\rr{n}} \bi{r}_{\rr{n}} + T_{\rr{t}} \bi{r}_{\rr{t}} . 
  \label{eq:equilibrium.eq.3}
\end{align}

\subsection{Incremental problem as nonlinear cone complementarity problem}

We have seen that the Coulomb friction law and the unilateral contact 
condition are written as \eqref{eq:complementarity.form.1}, 
where $\hat{g}_{j}(\Delta\bi{u})$ and $\Delta\bi{u}_{\rr{t}j}$ are 
related to $\Delta\bi{u}$ by 
\eqref{eq:kinematics.1} and \eqref{eq:kinematics.2}. 
Also, the equilibrium equation is written as \eqref{eq:equilibrium.eq.3}. 

Consequently, the incremental problem can be formulated as follows: 
\begin{subequations}\label{eq:nonlinear.complementarity}%
  \begin{align}
    & K \, \Delta\bi{u} - \bi{p}
    = \sum_{j=1}^{n} \bi{t}_{\rr{n}j} r_{\rr{n}j} 
    + \sum_{j=1}^{n} T_{\rr{t}j} \bi{r}_{\rr{t}j} , 
    \label{eq:nonlinear.complementarity.1} \\
    & F^{*} \ni 
    \begin{bmatrix}
      -g_{j} - \mu \| T_{\rr{t}j}^{\top} \, \Delta\bi{u} \| 
      + \bi{t}_{\rr{n}j}^{\top} \, \Delta\bi{u} \\
      T_{\rr{t}j}^{\top} \, \Delta\bi{u} \\
    \end{bmatrix}
    \perp
    \begin{bmatrix}
      r_{\rr{n}j} \\ \bi{r}_{\rr{t}j} \\
    \end{bmatrix}
    \in F , 
    \quad j=1,\dots,c. 
    \label{eq:nonlinear.complementarity.2}
  \end{align}
\end{subequations}
Here, $\Delta\bi{u}$, $r_{\rr{n}j}$, and $\bi{r}_{\rr{t}j}$ 
$(j=1,\dots,c)$ are unknown variables. 

It is worth noting that \eqref{eq:nonlinear.complementarity} is a 
nonlinear cone complementarity problem. 
Also, it is known that there exists no optimization problem the 
optimality condition of which corresponds to 
\eqref{eq:nonlinear.complementarity} (this is because, for the Coulomb 
friction law, there does not exist 
a {\em potential\/} \citep[section~3.9.2]{AB08}).

\section{Accelerated primal-dual algorithm}
\label{sec:algorithm}

In this section, based on the primal-dual algorithm \citep{CP16b} we 
develop an algorithm for solving problem \eqref{eq:nonlinear.complementarity}. 

\subsection{Optimization problem associated with \eqref{eq:nonlinear.complementarity}}

Since directly solving problem \eqref{eq:nonlinear.complementarity} 
seems not to be easy, we first consider an optimization problem that is 
similar to the one found in \cite{ACLM11}. 
This optimization problem is suited for application of a primal-dual 
algorithm. 

Let $\tilde{g}_{j} \in \Re$ $(j=1,\dots,c)$ be constants. 
Replace $g_{j} + \mu \| T_{\rr{t}j}^{\top} \Delta\bi{u} \|$ in 
\eqref{eq:nonlinear.complementarity} with $\tilde{g}_{j}$ to obtain the 
following complementarity problem: 
\begin{subequations}\label{eq:modified.complementarity}%
  \begin{align}
    & K \, \Delta\bi{u} - \bi{p}
    = \sum_{j=1}^{c} \bi{t}_{\rr{n}j} r_{\rr{n}j} 
    + \sum_{j=1}^{c} T_{\rr{t}j} \bi{r}_{\rr{t}j} , 
    \label{eq:modified.complementarity.1} \\
    & F^{*} \ni 
    \begin{bmatrix}
      -\tilde{g}_{j} + \bi{t}_{\rr{n}j}^{\top} \, \Delta\bi{u} \\
      T_{\rr{t}j}^{\top} \, \Delta\bi{u} \\
    \end{bmatrix}
    \perp
    \begin{bmatrix}
      r_{\rr{n}j} \\ \bi{r}_{\rr{t}j} \\
    \end{bmatrix}
    \in F , 
    \quad j=1,\dots,c. 
    \label{eq:modified.complementarity.2}
  \end{align}
\end{subequations}
For notational simplicity, define $\pi : \Re^{d} \to \Re$ by 
\begin{align*}
  \pi(\bi{u}) 
  = \frac{1}{2} \bi{u}^{\top} K \bi{u} - \bi{p}^{\top} \bi{u} . 
\end{align*}
We readily see that \eqref{eq:modified.complementarity} corresponds to 
the optimality condition for the following convex optimization 
problem (see appendix~\ref{sec:optimality} for details): 
\begin{align}
  \MIN_{\Delta\bi{u}} \quad
  \pi(\Delta\bi{u}) 
  + \sum_{j=1}^{c} \delta_{F^{*}}
  (-\tilde{g}_{j} + \bi{t}_{\rr{n}j}^{\top} \, \Delta\bi{u}, 
  T_{\rr{t}j}^{\top} \, \Delta\bi{u})  . 
  \label{P.approximation.3}
\end{align}
Here, for notational simplicity, 
for $v_{\rr{n}} \in \Re$, 
$\bi{v}_{\rr{t}} \in \Re^{2}$, 
and $\delta_{F^{*}} : \Re^{3} \to \overRe$ 
we write $\delta_{F^{*}}(v_{\rr{n}},\bi{v}_{\rr{t}})$ 
instead of 
$\delta_{F^{*}}\bigl( (v_{\rr{n}},\bi{v}_{\rr{t}}^{\top})^{\top} \bigr)$. 
Problem~\eqref{P.approximation.3} is a minimization problem of a convex 
quadratic function under second-order cone constraints. 

It follows from \eqref{eq:dual.cone.1} that we have 
\begin{align}
  \delta_{F^{*}}(\bi{v}) 
  &= -\inf_{\bi{r}_{j} \in \Re^{3}} 
  \{ \langle \bi{v},\bi{r}_{j} \rangle  + \delta_{F}(\bi{r}_{j}) \}  
  \notag\\
  &= \sup_{\bi{r}_{j} \in \Re^{3}} 
  \{ -\langle \bi{v},\bi{r}_{j} \rangle - \delta_{F}(\bi{r}_{j}) \} .  
  \label{eq:dual.cone.2}
\end{align}
Application of \eqref{eq:dual.cone.2} reduces \eqref{P.approximation.3} 
to the following form: 
\begin{align}
  \MIN_{\Delta\bi{u} \in \Re^{d}} \quad
  \pi(\Delta\bi{u}) + \sup_{\bi{r} \in \Re^{3c}} \left\{
  - \sum_{j=1}^{c} \left\langle
  \begin{bmatrix}
    r_{\rr{n}j} \\
    \bi{r}_{\rr{t}j} \\
  \end{bmatrix}
  , 
  \begin{bmatrix}
    -\tilde{g}_{j} 
    + \bi{t}_{\rr{n}j}^{\top} \, \Delta\bi{u} \\
    T_{\rr{t}j}^{\top} \, \Delta\bi{u} 
  \end{bmatrix}
  \right\rangle
  - \sum_{j=1}^{c} \delta_{F}(r_{\rr{n}j}, \bi{r}_{\rr{t}j}) 
  \right\} . 
  \label{P.approximation.4}
\end{align}

It should be clear that problems 
\eqref{eq:modified.complementarity}, 
\eqref{P.approximation.3}, and 
\eqref{P.approximation.4} in this section are equivalent to each other, 
but they are different from \eqref{eq:nonlinear.complementarity}.

\subsection{Primal-dual algorithm for optimization problem \eqref{P.approximation.4}}
\label{sec:algorithm.optimization}

For ease of comprehension of the algorithm presented in 
section~\ref{sec:algorithm.friction}, in this section we apply a 
primal-dual algorithm to problem \eqref{P.approximation.4}. 
Specifically, we apply \citet[Algorithm~1]{CP16b} 
(see also \citet[Algorithm~7]{CP16}). 

For notational simplicity, define $T \in \Re^{d \times 3c}$, 
$f : \Re^{3c} \to \overRe$, and $h : \Re^{3c} \to \Re$ by 
\begin{align*}
  T &= 
  \begin{bmatrix}
    \bi{t}_{\rr{n}1} & T_{\rr{t}1} & \cdots 
    & \bi{t}_{\rr{n}c} & T_{\rr{t}c} \\
  \end{bmatrix}
  , \\
  f^{*}(\bi{r}) 
  &= \sum_{j=1}^{c} 
  \delta_{F}(r_{\rr{n}j},\bi{r}_{\rr{t}j}) , \\
  h(\bi{r}) 
  &= \tilde{\bi{g}}^{\top} \bi{r}_{\rr{n}} . 
\end{align*}
Problem \eqref{P.approximation.4} is concisely written as follows: 
\begin{align}
  \min_{\Delta\bi{u}} \max_{\bi{r}} 
  \Bigl\{
  -\langle T^{\top} \, \Delta\bi{u}, \bi{r} \rangle 
  + \pi(\Delta\bi{u}) + h(\bi{r}) - f^{*}(\bi{r}) 
  \Bigr\} . 
  \label{P.approximation.6}
\end{align}
The primal-dual algorithm solving this problem updates the 
incumbent solution, denoted by $\bi{u}^{(k)}$ and $\bi{r}^{(k)}$, as 
\begin{align}
  \bi{r}^{(k+1)} 
  &:= \prox_{\alpha f^{*}}
  (\bi{r}^{(k)} 
  + \alpha (\nabla h(\bi{r}^{(k)}) 
  - T^{\top} \, \Delta\hat{\bi{u}}^{(k)})) , 
  \label{p.d.friction.1} \\
  \Delta\bi{u}^{(k+1)} 
  &:= \prox_{\beta \pi}
  (\Delta\bi{u}^{(k)} + \beta T \bi{r}^{(k+1)}) ,  
  \label{p.d.friction.2} \\
  \Delta\hat{\bi{u}}^{(k+1)} 
  & := \Delta\bi{u}^{(k+1)} 
  + \theta (\Delta\bi{u}^{(k+1)} - \Delta\bi{u}^{(k)}) . 
  \label{p.d.friction.3}
\end{align}
Here, $\alpha>0$ and $\beta>0$ are step lengths, 
and $\theta \in [0,1]$ is a constant. 

A direct calculation using definition \eqref{eq:def.proximal.mapping} 
of the proximal mapping yields 
\begin{align*}
  \prox_{\alpha f^{*}}(\bi{r}) &= 
  \begin{bmatrix}
    \Pi_{F}(\bi{r}_{1}) \\
    \vdots \\
    \Pi_{F}(\bi{r}_{c}) \\
  \end{bmatrix}
  , \\
  \prox_{\beta \pi}(\bi{u}) 
  &= (\beta K + I)^{-1} (\bi{u} + \beta \bi{p}) . 
\end{align*}
Therefore, the updates in \eqref{p.d.friction.1}, \eqref{p.d.friction.2}, 
and \eqref{p.d.friction.3} can be described explicitly as 
\refalg{alg:frictional.2}. 

\begin{algorithm}
  \caption{Primal-dual algorithm for optimization problem 
  \eqref{eq:modified.complementarity}. }
  \label{alg:frictional.2}
  \begin{algorithmic}[1]
    \Require
    $\Delta\bi{u}^{(0)} \in \Re^{d}$, 
    $\bi{r}^{(0)} \in \Re^{3c}$, 
    $\alpha >0$, $\beta > 0$, $\theta \in [0,1]$. 
    \State
    $\Delta\hat{\bi{u}}^{(0)} \gets \Delta\bi{u}^{(0)}$. 
    \For{$k=0,1,2,\dots$}
    \State
    $\bi{s}_{\rr{n}}  \gets 
    \bi{r}_{\rr{n}}^{(k)} 
    + \alpha (\tilde{\bi{g}} - T_{\rr{n}}^{\top} \, \Delta\hat{\bi{u}}^{(k)})$. 
    \State
    $\bi{s}_{\rr{t}}  \gets 
    \bi{r}_{\rr{t}}^{(k)} - \alpha T_{\rr{t}}^{\top} \, \Delta\hat{\bi{u}}^{(k)}$. 
    \State
    $\bi{r}_{j}^{(k+1)}  \gets \Pi_{F}(\bi{s}_{j})$  
    $(j=1,\dots,c)$. 
    \State
    $\bi{b}  \gets 
    \Delta\bi{u}^{(k)} + \beta (T \bi{r}^{(k+1)} + \bi{p} )$. 
    \State
    Solve 
    $(\beta K + I) \, \Delta\bi{u}^{(k+1)} = \bi{b}$ 
    to obtain $\Delta\bi{u}^{(k+1)}$. 
    \State
    $\Delta\hat{\bi{u}}^{(k+1)} \gets 
    \Delta\bi{u}^{(k+1)} + \theta(\Delta\bi{u}^{(k+1)} - \Delta\bi{u}^{(k)})$. 
    \EndFor
  \end{algorithmic}
\end{algorithm}

\begin{remark}
  Problem \eqref{P.approximation.3} is regarded as a particular case of 
  the following convex optimization problem in variable $\bi{x} \in \Re^{n}$: 
  \begin{subequations}\label{P.generalized.problem.1}%
    \begin{alignat}{3}
      & \MIN
      &{\quad}& 
      \hat{f}(\bi{x}) 
      := \frac{1}{2} \bi{x}^{\top} Q \bi{x} + \bi{h}^{\top} \bi{x}  \\
      & \ST && 
      A_{i} \bi{x} + \bi{b}_{i} \in L^{n_{i}} , 
      \quad i=1,\dots,m, \\
      & && 
      C \bi{x} + \bi{d} = \bi{0}. 
    \end{alignat}
  \end{subequations}
  Here, $Q \in \Re^{n \times n}$ is symmetric and positive semidefinite, and 
  $L^{n_{i}}$ denotes the $n_{i}$-dimensional second-order cone, i.e., 
  \begin{align*}
    L^{n_{i}} = \{
    (x_{0},\bi{x}_{1}) \in \Re \times \Re^{n_{i}-1} 
    \mid
    x_{0} \ge \| \bi{x}_{1} \| 
    \} . 
  \end{align*}
  Problem \eqref{P.generalized.problem.1} is minimization of a convex 
  quadratic function under second-order cone constraints. 
  It follows from the self-duality of the second-order cone that 
  problem \eqref{P.generalized.problem.1} is equivalently rewritten as 
  follows: 
  \begin{align}
    \MIN_{\bi{x}} 
    \quad 
    \sup_{\bi{s}_{1} \in L^{n_{1}},\dots,\bi{s}_{m} \in L^{n_{m}},\bi{y}}
    \Bigl\{
    \frac{1}{2} \bi{x}^{\top} Q \bi{x} + \bi{p}^{\top} \bi{x}  
    - \sum_{i=1}^{m} 
    \langle \bi{s}_{i}, A_{i} \bi{x} + \bi{b}_{i} \rangle
    - \bi{y}^{\top} (C \bi{x} + \bi{d})
    \Bigr\} . 
    \label{P.generalized.problem.min.max}
  \end{align}
  The primal-dual algorithm for solving this problem 
  updates the dual variables as 
  \begin{align}
    \bi{s}^{(k+1)} 
    & := \Pi_{L^{n_{i}}} 
    (\bi{s}_{i}^{(k)} - \alpha (A_{i} \bi{x}^{(k)} + \bi{b}_{i})) ,
    \quad i=1,\dots,m, 
    \label{eq:generalized.problem.2} \\
    \bi{y}^{(k+1)}
    & := \bi{y}^{(k)} - \alpha (C \bi{x}^{(k)} + \bi{d}) , 
    \label{eq:generalized.problem.3}
  \end{align}
  and then updates the primal variable as 
  \begin{align*}
    \bi{x}^{(k+1)} 
    &:= \prox_{\beta \hat{f}}
    \Bigl( \bi{x}^{(k)} + \sum_{i=1}^{n} A_{i}^{\top} \bi{s}_{i}^{(k+1)}
    + C^{\top} \bi{y}^{(k+1)} \Bigr) , 
  \end{align*}
  where 
  \begin{align*}
    \prox_{\beta \hat{f}}(\bi{x}) 
    = (\beta Q + I)^{-1} (\bi{x} + \beta \bi{h}) . 
  \end{align*}
  This algorithm converges to a saddle point of problem 
  \eqref{P.generalized.problem.min.max} \citep{CP16b}. 
  An advantage of using the primal-dual algorithm lies in the fact that 
  the update of the dual variables in \eqref{eq:generalized.problem.2} 
  and \eqref{eq:generalized.problem.3} can be performed very easily 
  (particularly an explicit formula for the projection in 
  \eqref{eq:generalized.problem.2} is available \citep{FLT01}), 
  compared with the projection of the primal variable $\bi{x}$ onto the 
  feasible set of problem \eqref{P.generalized.problem.1}. 
  \refalg{alg:frictional.2} presented in this section has the same 
  advantage, compared with directly handling 
  the primal formulation in \eqref{P.approximation.3}.  
  \finbox
\end{remark}

\subsection{Algorithm for frictional contact problem}
\label{sec:algorithm.friction}

We have seen in section~\ref{sec:algorithm.optimization} that problem 
\eqref{eq:modified.complementarity} can be solved with 
\refalg{alg:frictional.2}. 
To deal with the frictional contact problem in 
\eqref{eq:nonlinear.complementarity}, 
we make two alterations to \refalg{alg:frictional.2} as follows. 

One alteration is to implement an acceleration scheme. 
Since $\pi$ is a strongly convex function, the acceleration scheme in 
\citet[Algorithm~4]{CP16b} 
(see also \citet[Algorithm~8]{CP16}) 
is likely to work efficient. 


The other is to update $\tilde{g}_{j}$ $(j=1,\dots,c)$ 
at each iteration. 
Since we have obtained problem \eqref{eq:modified.complementarity} by 
replacing $g_{j} + \mu \| T_{\rr{t}j}^{\top} \, \Delta\bi{u} \|$ in 
\eqref{eq:nonlinear.complementarity} with $\tilde{g}_{j}$, 
a natural update may be 
\begin{align*}
  \tilde{g}_{j} 
  := g_{j} + \mu \| T_{\rr{t}j}^{\top} \, \Delta\bi{u}^{(k)} \| , 
  \quad j=1,\dots,c . 
\end{align*}
It is worth noting that, with this update, there is no guarantee of 
convergence to a solution of problem \eqref{eq:nonlinear.complementarity}. 
In the numerical experiments reported in section~\ref{sec:ex}, 
the proposed algorithm converges to a solution of every problem 
instance. 

Application of the two alterations above to \refalg{alg:frictional.2} 
yields \refalg{alg:frictional.accelerated.2}. 
Here, $\sigma_{T}$ and $\mu_{\pi}$ are the maximum singular value of $T$ 
and the minimum eigenvalue of $\nabla^{2}\pi(\bi{x}) = K$, respectively. 

In line~\ref{alg:frictional.accelerated.equation} 
of \refalg{alg:frictional.accelerated.2}, we solve a system of linear 
equations to obtain $\Delta\bi{u}^{(k+1)}$. 
We adopt a preconditioned conjugate gradient method 
with setting $\Delta\bi{u}^{(k)}$ 
as an initial point, because we may probably expect that change 
from $\Delta\bi{u}^{(k)}$ to $\Delta\bi{u}^{(k+1)}$ is not large. 
The computation in line~\ref{alg:frictional.accelerated.2.projection} is 
performed for each contact candidate node $j=1,\dots,c$ independently, 
and hence can be highly parallelized. 
In section~\ref{sec:algorithm.projection}, we present 
a concrete procedure for computing $\Pi_{F}(\bi{s}_{j})$. 
All the other computations in \refalg{alg:frictional.accelerated.2} are 
additions and multiplications, which are computationally cheap. 

%

\begin{algorithm}
  \caption{Primal-dual algorithm for frictional contact problem 
  \eqref{eq:nonlinear.complementarity}. }
  \label{alg:frictional.accelerated.2}
  \begin{algorithmic}[1]
    \Require
    $\Delta\bi{u}^{(0)}$, 
    $\bi{r}^{(0)}$, $\alpha_{0} > 0$. 
    \State
    $\Delta\hat{\bi{u}}^{(0)} \gets \Delta\bi{u}^{(0)}$. 
    \State
    $\beta_{0} \gets 1 / (\alpha_{0} \sigma_{T}^{2})$. 
    \For{$k=0,1,2,\dots$}
    \State
    $\tilde{g}_{j}^{(k)} 
    \gets g_{j} + \mu \| T_{\rr{t}j}^{\top} \, \Delta\bi{u}^{(k)} \|$ 
    $(j=1,\dots,c)$. 
    \State
    $\bi{s}_{\rr{n}}  \gets 
    \bi{r}_{\rr{n}}^{(k)} 
    + \alpha_{k} 
    (\tilde{\bi{g}}^{(k)} - T_{\rr{n}}^{\top} \, \Delta\hat{\bi{u}}^{(k)})$. 
    \State
    $\bi{s}_{\rr{t}}  \gets 
    \bi{r}_{\rr{t}}^{(k)} 
    - \alpha_{k} T_{\rr{t}}^{\top} \, \Delta\hat{\bi{u}}^{(k)}$. 
    \State \label{alg:frictional.accelerated.2.projection}
    $\bi{r}_{j}^{(k+1)}  \gets \Pi_{F}(\bi{s}_{j})$  
    $(j=1,\dots,c)$. 
    \State
    $\bi{b}  \gets 
    \Delta\bi{u}^{(k)} + \beta_{k} (T \bi{r}^{(k+1)} + \bi{p} )$. 
    \State \label{alg:frictional.accelerated.equation}
    Solve 
    $(\beta_{k} K + I) \, \Delta\bi{u}^{(k+1)} = \bi{b}$ 
    to obtain $\Delta\bi{u}^{(k+1)}$. 
    \State
    $\displaystyle \theta_{k} 
    \gets \frac{1}{\sqrt{1 + \mu_{\pi} \beta_{k}}}$. 
    \State
    $\displaystyle \alpha_{k+1} 
    \gets \frac{\alpha_{k}}{\theta_{k}}$, 
    $\beta_{k+1} 
    \gets \theta_{k} \beta_{k}$. 
    \State
    $\Delta\hat{\bi{u}}^{(k+1)} \gets 
    \Delta\bi{u}^{(k+1)} 
    + \theta_{k} (\Delta\bi{u}^{(k+1)} - \Delta\bi{u}^{(k)})$. 
    \EndFor
  \end{algorithmic}
\end{algorithm}

\subsection{Projection onto Coulomb's friction cone}
\label{sec:algorithm.projection}

\begin{figure*}[tbp]
  \centering
  \includegraphics[scale=0.50]{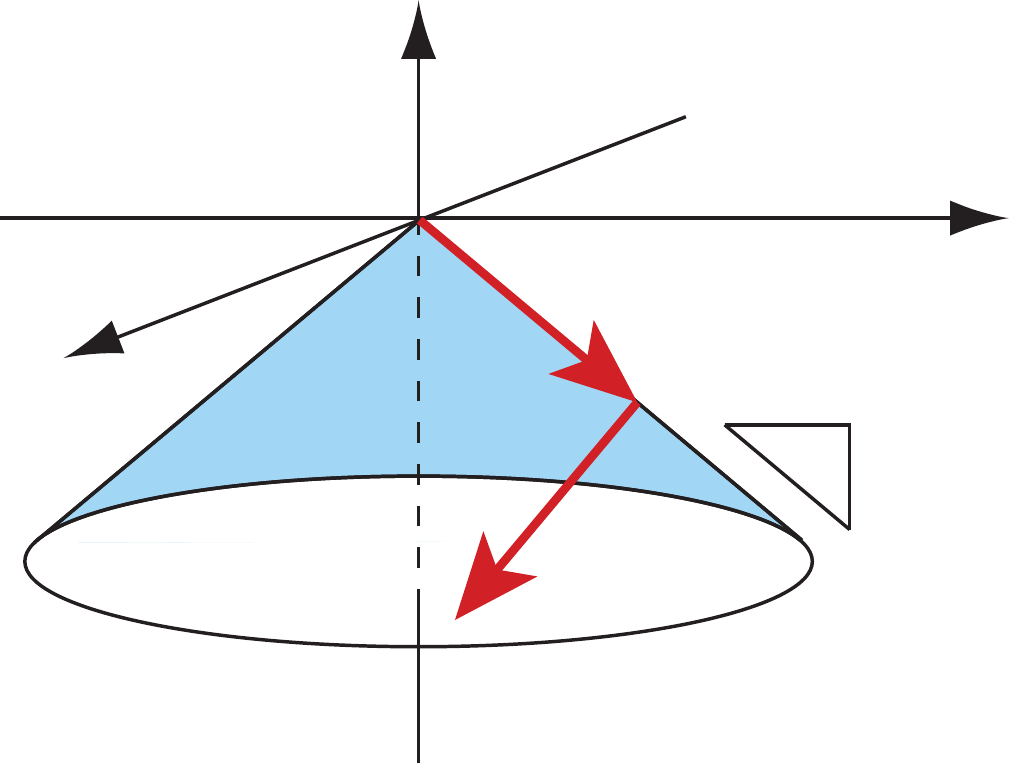}
  \begin{picture}(0,0)
    \put(-200,-10){
    \put(98,113.5){\footnotesize $r_{\rr{n}}$ }
    \put(47.5,66){\footnotesize $r_{\rr{t}1}$ }
    \put(197,86){\footnotesize $r_{\rr{t}2}$ }
    \put(140,70){\footnotesize $\bi{\xi}_{2}$ }
    \put(128,33){\footnotesize $\bi{\xi}_{1}$ }
    \put(163,61){\footnotesize $1$ }
    \put(175,50){\footnotesize $\mu$ }
    \put(49,43){\footnotesize $F$ }
    }
  \end{picture}
  \caption{Projection onto the Coulomb friction cone.}
  \label{fig:friction_projection}
\end{figure*}

In line~\ref{alg:frictional.accelerated.2.projection} of 
\refalg{alg:frictional.accelerated.2}, we compute the projection of 
$\bi{s}_{j} \in \Re^{3}$ onto the Coulomb friction cone $F$. 
This can be easily performed as follows. 

For notational simplicity, we consider computation of 
$\Pi_{F}(s_{\rr{n}},\bi{s}_{\rr{t}})$, i.e., we omit subscript $j$. 
Define $\bi{\xi}_{1}$, $\bi{\xi}_{2} \in \Re^{3}$ by 
\begin{align*}
  \bi{\xi}_{1} = -\frac{1}{1+\mu^{2}} 
  \begin{bmatrix}
    \mu \\[1.0ex]
    \displaystyle
    \frac{\bi{s}_{\rr{t}}}{\| \bi{s}_{\rr{t}} \|} \\
  \end{bmatrix}
  , \quad
  \bi{\xi}_{2} = \frac{1}{1+\mu^{2}} 
  \begin{bmatrix}
    -1 \\[1.0ex]
    \displaystyle
    \mu \frac{\bi{s}_{\rr{t}}}{\| \bi{s}_{\rr{t}} \|} \\
  \end{bmatrix}
  ; 
\end{align*}
see \reffig{fig:friction_projection}.\footnote{%
When $\bi{s}_{\rr{t}}=\bi{0}$, we use any unit vector instead of 
$\bi{s}_{\rr{t}}/\| \bi{s}_{\rr{t}} \|$. 
It is worth noting that for implementation we do not need to consider 
this case; see \refalg{alg:projection.Coulomb}. } 
Also, define $\lambda_{1}$, $\lambda_{2} \in \Re$ by 
\begin{align*}
  \lambda_{1} = -\mu s_{\rr{n}} - \| \bi{s}_{\rr{t}} \| , 
  \quad
  \lambda_{2} = -s_{\rr{n}} + \mu \| \bi{s}_{\rr{t}} \| . 
\end{align*}
Then we have 
\begin{align}
  \Pi_{F}(s_{\rr{n}}, \bi{s}_{\rr{t}}) 
  = \max \{ 0, \lambda_{1} \} \bi{\xi}_{1} 
  + \max \{ 0, \lambda_{2} \} \bi{\xi}_{2} . 
  \label{eq:projection.friction.cone.1}
\end{align}
It is worth noting that, when $\mu=1$, 
\eqref{eq:projection.friction.cone.1} corresponds to a formula of the 
projection onto the second-order cone \citep{FLT01} based on the 
spectral factorization in the Jordan algebra. 

The computational procedure of the projection is described in 
\refalg{alg:projection.Coulomb}. 

\begin{algorithm}
  \caption{Projection onto the Coulomb's friction cone, $\Pi_{F}(s_{\rr{n}},\bi{s}_{\rr{t}})$. }
  \label{alg:projection.Coulomb}
  \begin{algorithmic}[1]
    \Require
    $(s_{\rr{n}},\bi{s}_{\rr{t}}) \in \Re \times \Re^{2}$. 
    \State
    $\lambda_{1} \gets -\mu s_{\rr{n}} - \| \bi{s}_{\rr{t}} \|$. 
    \If{$\lambda_{1} \ge 0$}
    \State
    $\Pi_{F}(s_{\rr{n}},\bi{s}_{\rr{t}}) 
    \gets  (s_{\rr{n}},\bi{s}_{\rr{t}})$. 
    \Else
    \State
    $\lambda_{2} \gets -s_{\rr{n}} + \mu \| \bi{s}_{\rr{t}} \|$. 
    \If{ $\lambda_{2} \le 0$ }
    \State
    $\Pi_{F}(s_{\rr{n}},\bi{s}_{\rr{t}}) \gets  (0,\bi{0})$. 
    \Else
    \State
    $\displaystyle
    \Pi_{F}(s_{\rr{n}},\bi{s}_{\rr{t}})  \gets
    \frac{\lambda_{2}}{1+\mu^{2}} \Bigl(
    -1 , \mu \frac{\bi{s}_{\rr{t}}}{\| \bi{s}_{\rr{t}} \|}
    \Bigr)$. 
    \EndIf
    \EndIf
  \end{algorithmic}
\end{algorithm}

\section{Numerical experiments}
\label{sec:ex}

In this section, we demonstrate efficiency of the proposed method 
through numerical experiments. 
Section~\ref{sec:ex.implementation} describes details of implementation. 
Sections~\ref{sec:ex.2D} and \ref{sec:ex.3D} report on two numerical 
examples. 
Computation was carried out on 
a 2.6{\,}GHz Intel Core i7-9750H processor with 32{\,}GB RAM. 
In the following, we omit units of physical quantities for simplicity. 
Young's modulus and Poisson's ratio of elastic bodies are 
$1$ and $0.3$, respectively.

\subsection{Implementation}
\label{sec:ex.implementation}

\refalg{alg:frictional.accelerated.2} was implemented 
in Matlab ver.~9.8.0. 
Initial values were set to 
$\Delta\bi{u}^{(0)}=\bi{0}$, $\bi{r}^{(0)} = \bi{0}$, and 
$\alpha_{0} = 10^{-1}$. 
Stopping criterion was 
$\| \Delta\bi{u}^{(k+1)} - \Delta\bi{u}^{(k)} \| \le \epsilon$ 
with $\epsilon=10^{-12}$. 

The projection onto the friction cone in 
line~\ref{alg:frictional.accelerated.2.projection} is computed with 
\refalg{alg:projection.Coulomb}. 
Matlab implementation of \refalg{alg:projection.Coulomb} was compiled 
into a MEX-file by using Matlab Coder, where the loop for $j=1,\dots,c$ 
was implemented with Matlab built-in function \texttt{parfor}. 
For this parallel loop computation, Matlab was allowed to use up to 6 cores. 

As explained in section~\ref{sec:algorithm.friction}, 
a preconditioned conjugate gradient method was used for solution of a system of 
linear equations in line~\ref{alg:frictional.accelerated.equation} 
of \refalg{alg:frictional.accelerated.2}. 
Matlab built-in function \texttt{pcg} was used, where the maximum number 
of iterations and the tolerance for termination were 
set to $10^{4}$ and $10^{-10}$, respectively. 
The minimum eigenvalue of $K$, i.e., $\mu_{\pi}$, was computed by 
\texttt{eigs($\,\cdot\,$,1,'smallestabs')}. 
The maximum singular value of $T$, i.e., $\sigma_{T}$, was computed by 
\texttt{svds($\,\cdot\,$,1,'largest')}, 
with setting the maximum number of iterations to $10^{8}$. 

For comparison, the same problem instances were solved with a 
regularized and smoothed Newton method \citep{HYF05}, which is 
implemented in ReSNA \cite{Hay20}. 
Specifically, we solved 
problem \eqref{P.linear.SOC.2} in appendix~\ref{sec:linear_form}, 
which is a second-order cone linear complementarity problem. 
Parameter \texttt{tole} of ReSNA was set to $10^{-7}$. 
Both for ReSNA and \refalg{alg:frictional.accelerated.2}, 
the stiffness matrix $K$ was stored as a Matlab sparse form. 

Besides computation time, we compare accuracy of the obtained solutions. 
With referring to problem \eqref{eq:nonlinear.complementarity}, 
we consider the following three residuals. 
First, the residual of \eqref{eq:nonlinear.complementarity.1} is 
\begin{align}
  \| K \Delta\bi{u}^{(k)} 
  - \bi{p} - T_{\rr{n}} \bi{r}_{\rr{n}}^{(k)} 
  - T_{\rr{t}j} \bi{r}_{\rr{t}}^{(k)} \| . 
  \label{eq:residual.1}
\end{align}
Second, as the residual of the complementarity conditions in 
\eqref{eq:nonlinear.complementarity.2}, we consider 
\begin{align}
  \left|
  \sum_{j=1}^{c}
  \left\langle
  \begin{bmatrix}
    -g_{j} - \mu \| T_{\rr{t}j}^{\top} \Delta\bi{u}^{(k)} \| 
    + \bi{t}_{\rr{n}j}^{\top} \Delta\bi{u}^{(k)} \\
    T_{\rr{t}j}^{\top} \Delta\bi{u}^{(k)} \\
  \end{bmatrix}
  ,
  \begin{bmatrix}
    r_{\rr{n}j}^{(k)} \\ \bi{r}_{\rr{t}j}^{(k)} \\
  \end{bmatrix}
  \right\rangle
  \right|
  = 
  | \langle \bi{u}^{(k)}, T \bi{r}^{(k)} \rangle 
  - \langle \Delta\tilde{\bi{g}}^{(k)}, \bi{r}_{\rr{n}}^{(k)} \rangle |
  . 
  \label{eq:residual.2}
\end{align}
Finally, for the inequality constraints 
$g_{j} - \bi{t}_{\rr{n}j}^{\top} \Delta\bi{u} \ge 0$ $(j=1,\dots,c)$ 
in \eqref{eq:nonlinear.complementarity.2}, the residual is 
\begin{align}
  \| \min \{ \bi{g} - T_{\rr{n}}^{\top} \Delta\bi{u}^{(k)} , \bi{0} \} \| . 
  \label{eq:residual.3}
\end{align}
It is worth noting that $\bi{r}_{j}^{(k)} \in F$ $(j=1,\dots,c)$ 
is satisfied at every 
iteration of \refalg{alg:frictional.accelerated.2}, due to the 
projection in line~\ref{alg:frictional.accelerated.2.projection}.

\subsection{Example (I): two-dimensional problem}
\label{sec:ex.2D}

\begin{figure*}[tbp]
  \centering
  \includegraphics[scale=0.50]{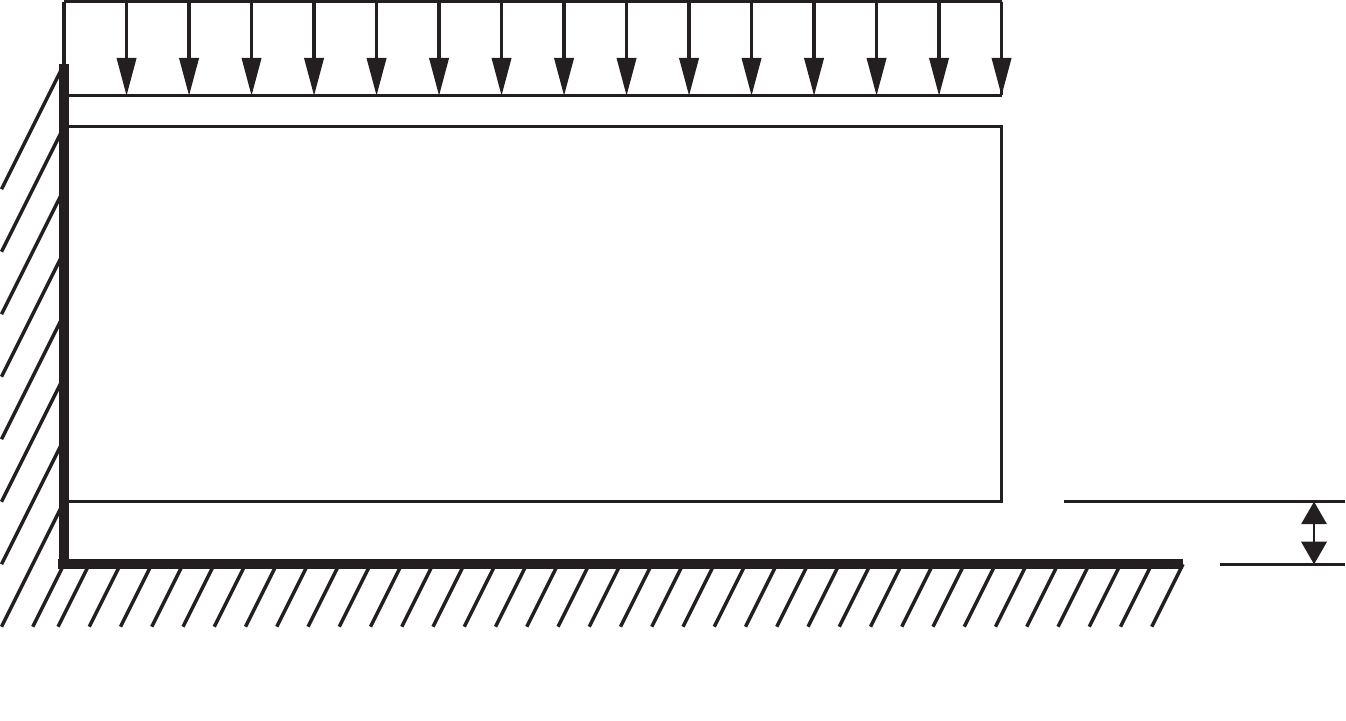}
  \begin{picture}(0,0)
    \put(-150,-150){
    \put(149,176){{\small $g_{j}$}}
    \put(30,152){{\small obstacle}}
    }
  \end{picture}
  \caption[]{Problem setting of example (I). }
  \label{fig:ex_body_2d}
\end{figure*}

\begin{figure*}[tbp]
  \centering
  \subfloat[]{
  \label{fig:time_log_dof_prog2ex_rev}
  \includegraphics[scale=0.50]{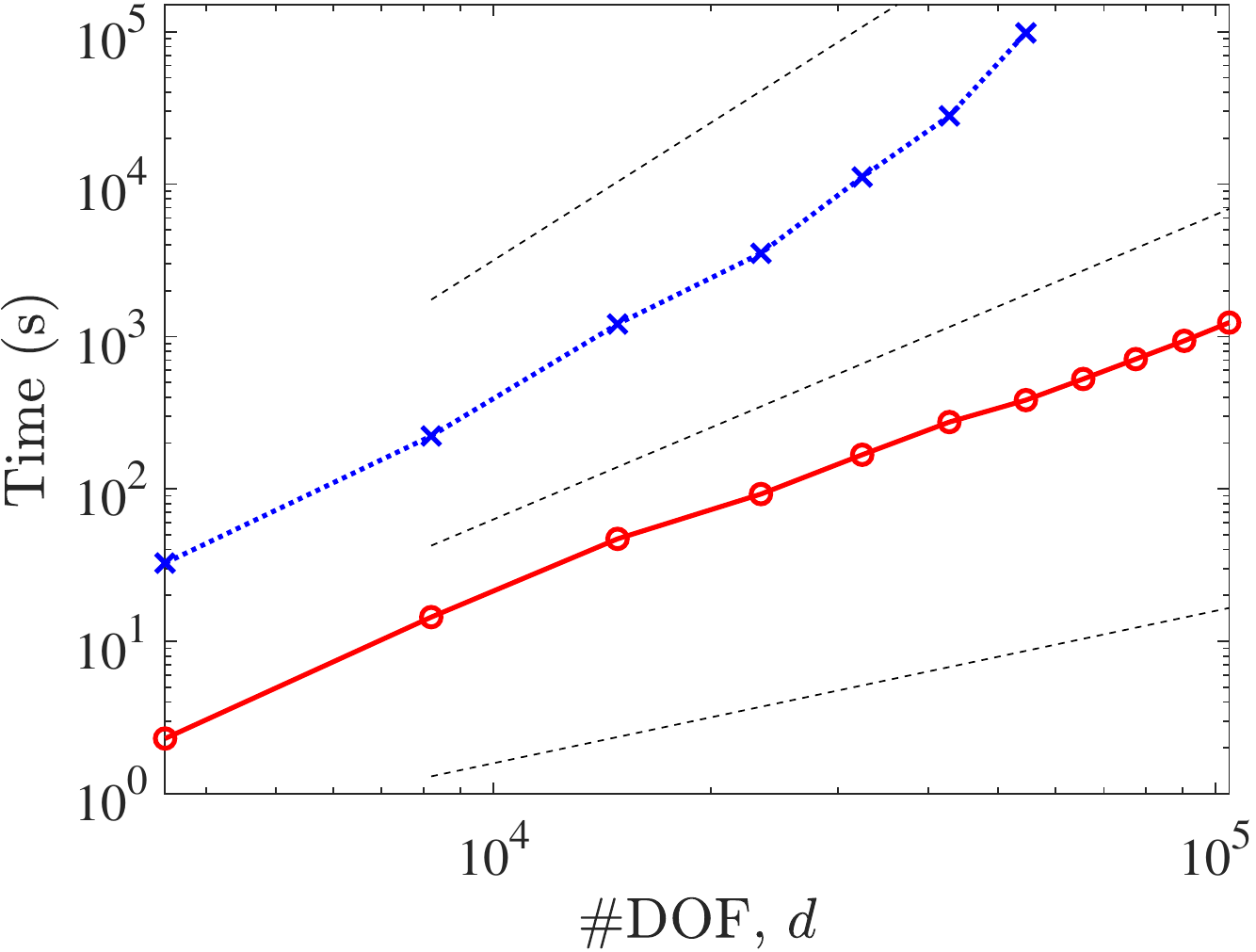}
  }
  \hfill
  \subfloat[]{
  \label{fig:iter_pd_prog2ex_rev}
  \includegraphics[scale=0.50]{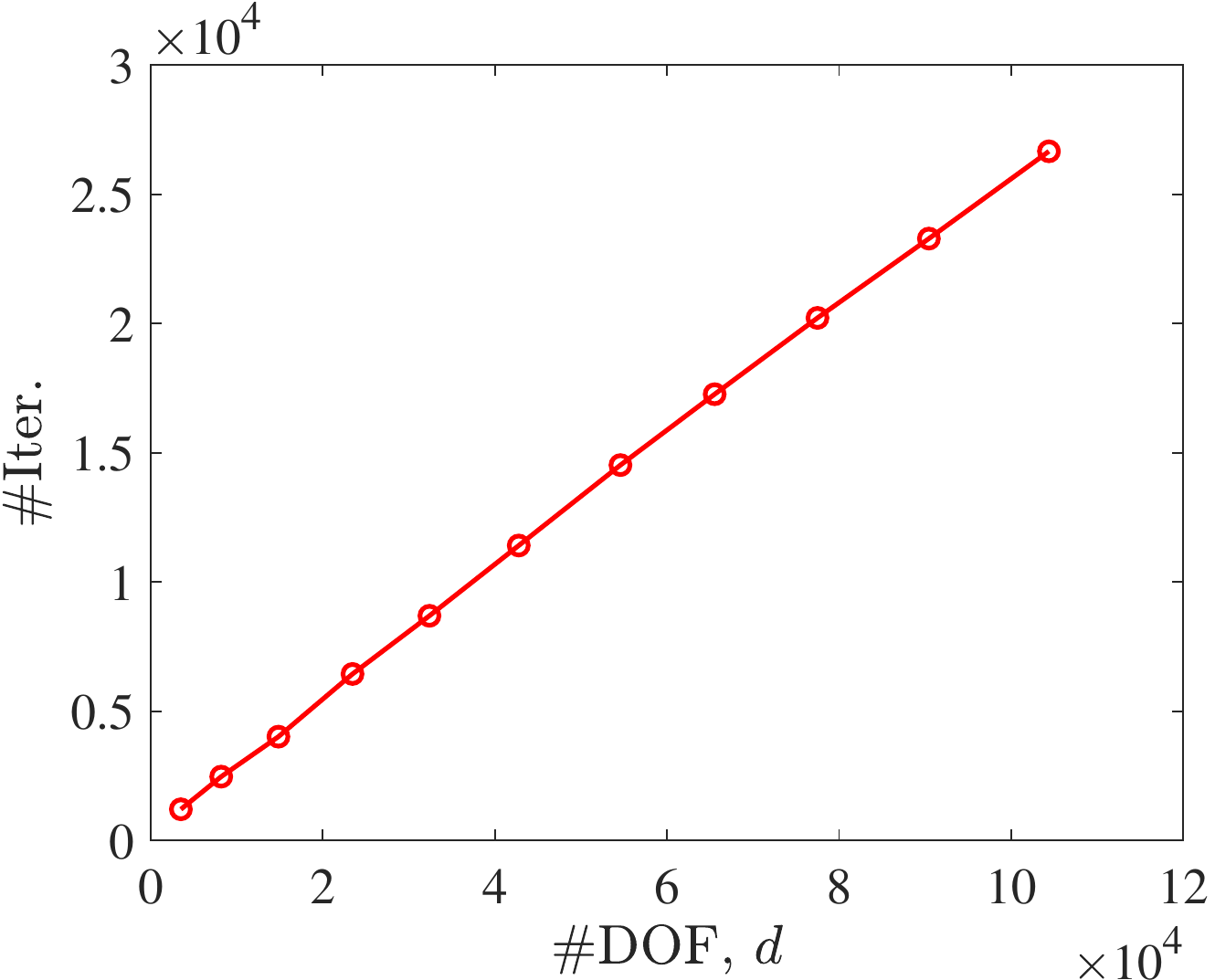}
  }
  \par
  \subfloat[]{
  \label{fig:iter_soc_prog2ex_rev}
  \includegraphics[scale=0.50]{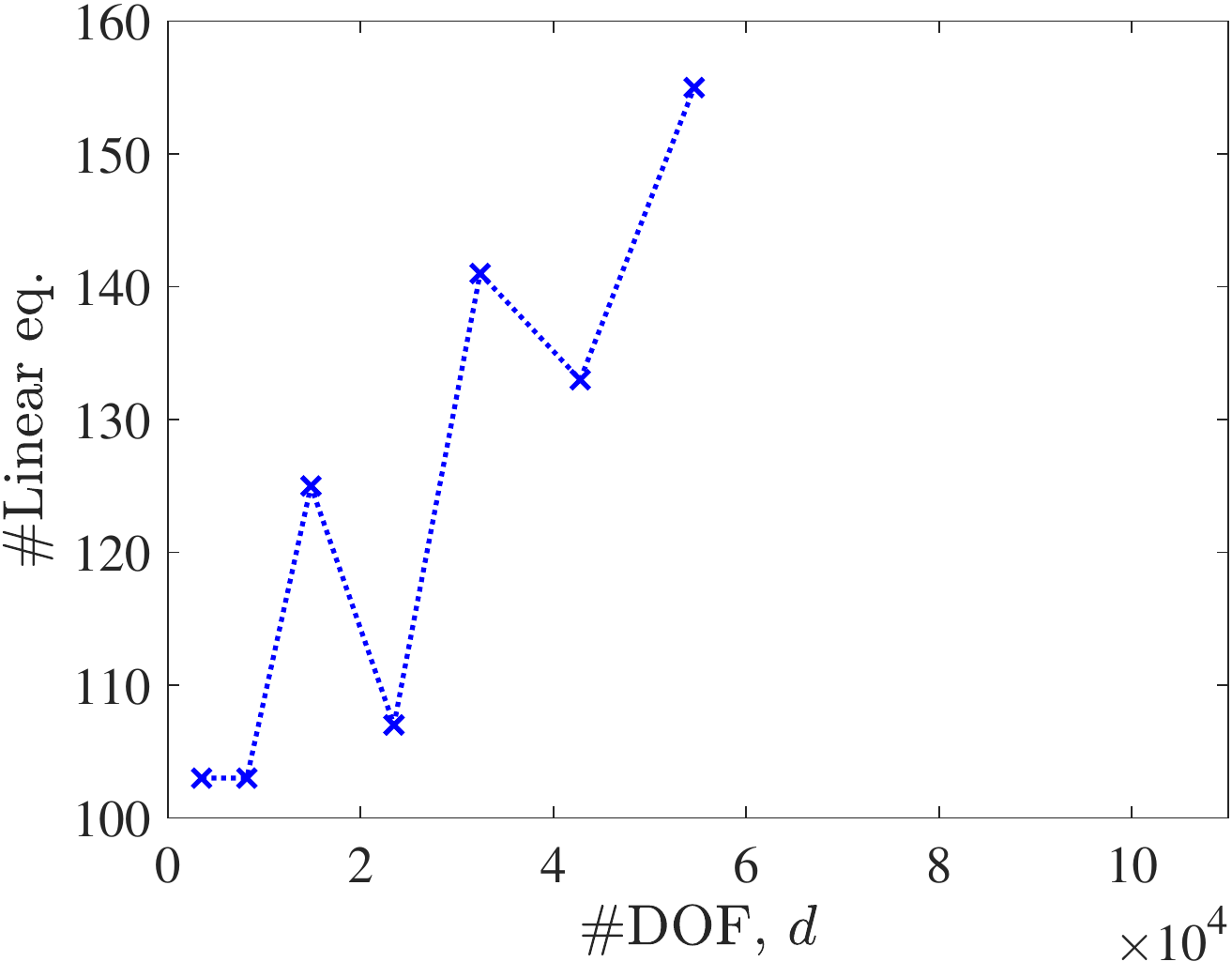}
  }
  \hfill
  \subfloat[]{
  \label{fig:resid_eq_prog2ex_rev}
  \includegraphics[scale=0.50]{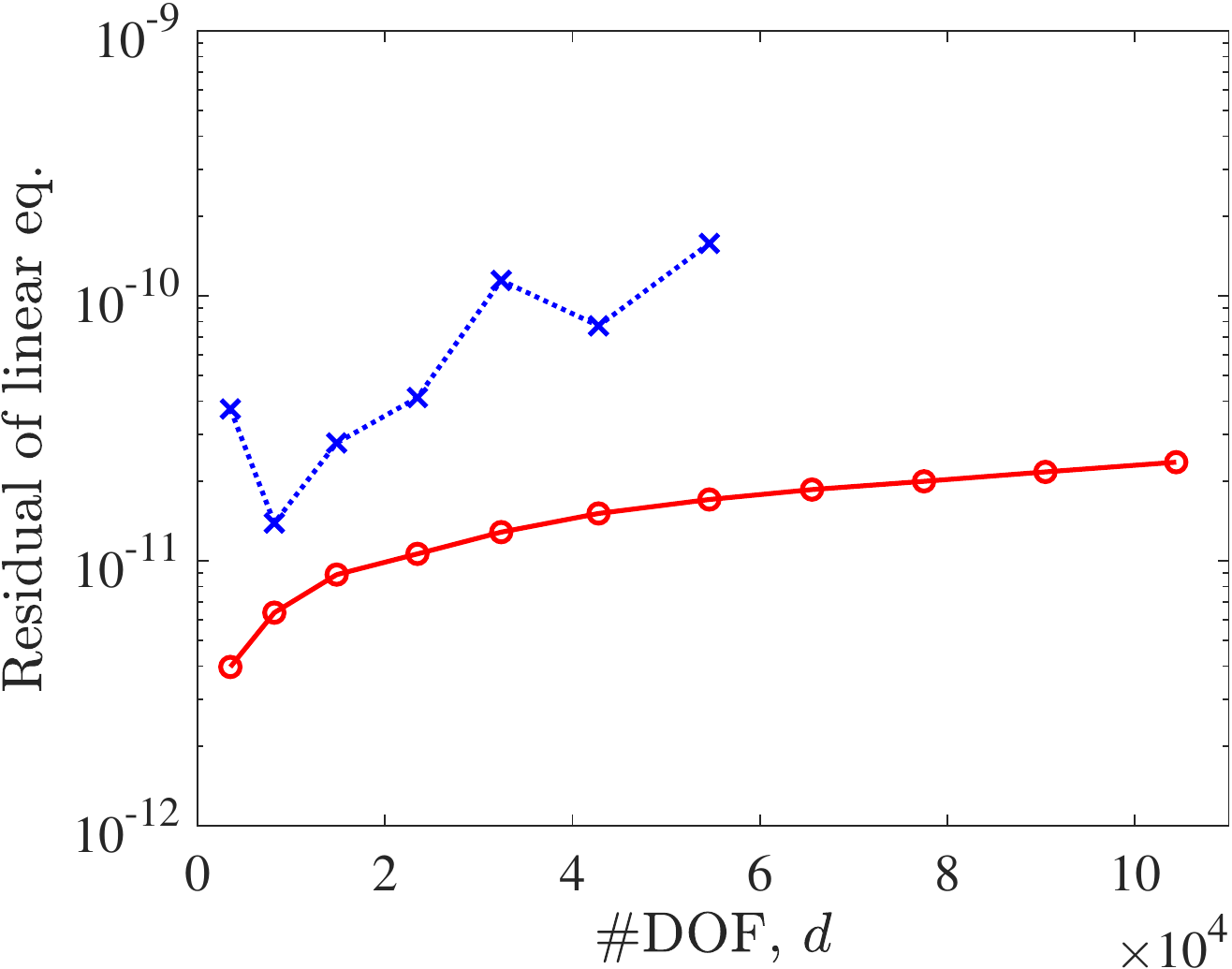}
  }
  \par
  \subfloat[]{
  \label{fig:resid_compl_prog2ex_rev}
  \includegraphics[scale=0.50]{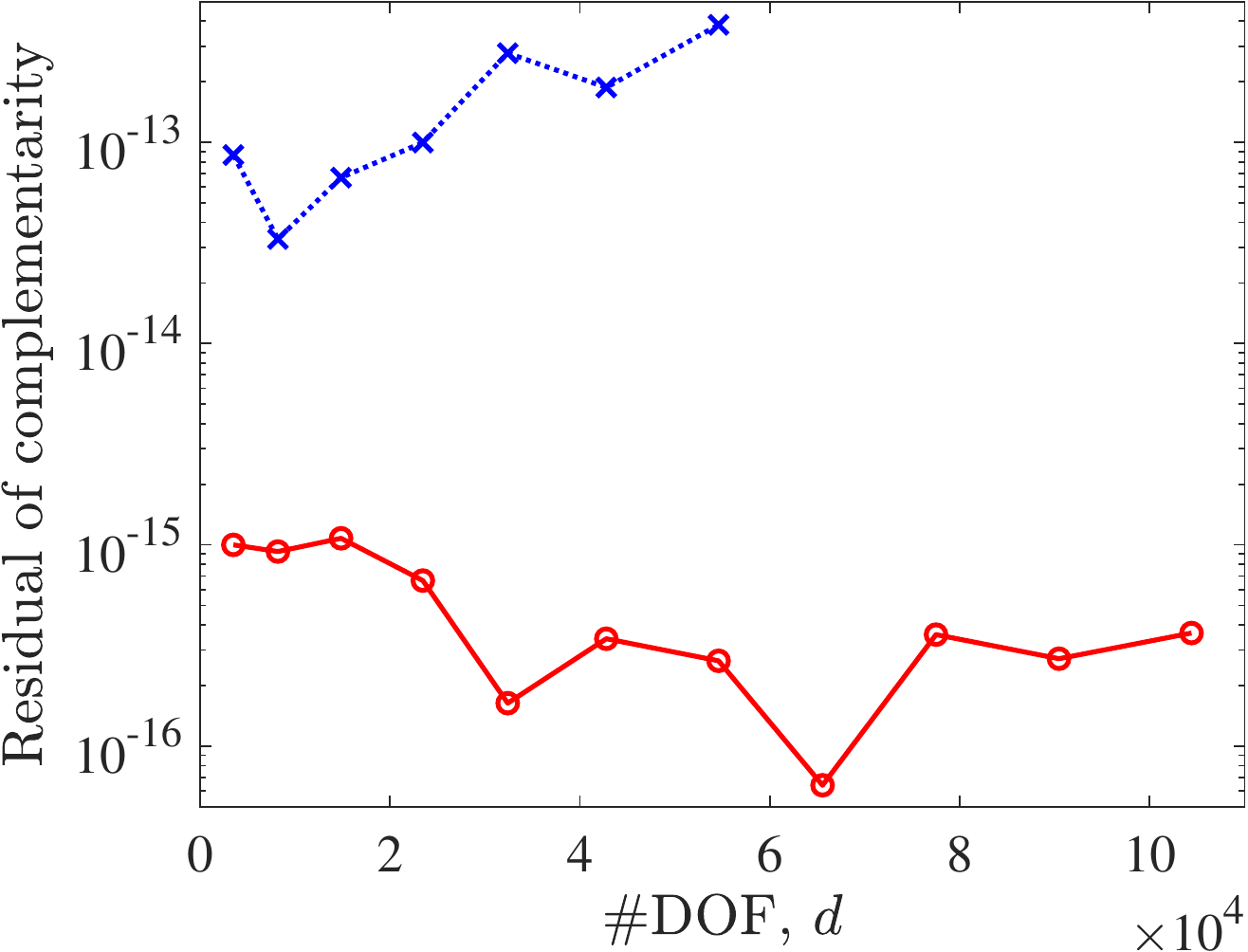}
  }
  \hfill
  \subfloat[]{
  \label{fig:resid_gap_ineq_prog2ex_rev}
  \includegraphics[scale=0.50]{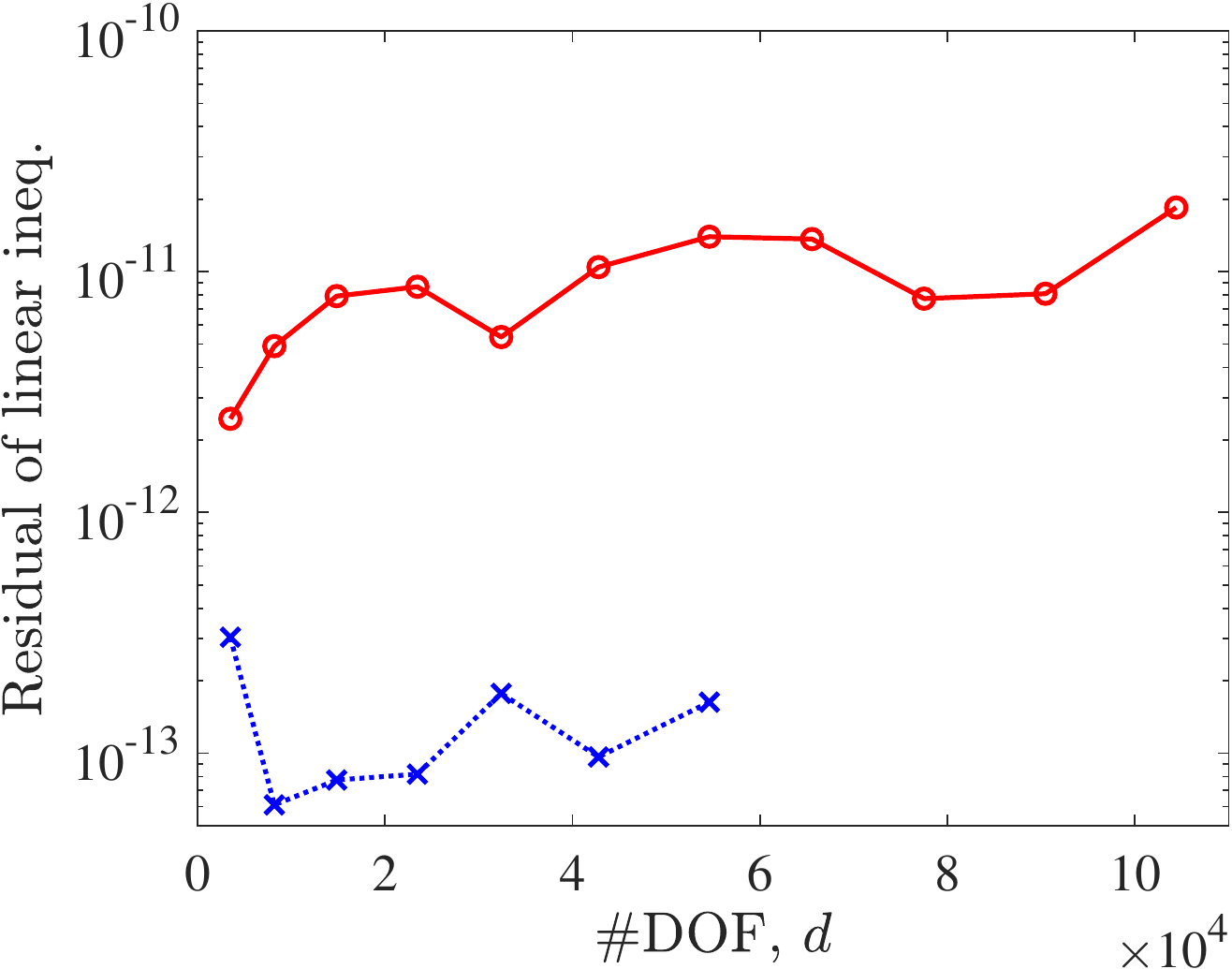}
  }
  \caption[]{Computational results of example (I). 
  ``{\em Solid line\/}'' The proposed method; 
  and ``{\em dotted line\/}'' ReSNA. 
  \subref{fig:time_log_dof_prog2ex_rev} The computation time; 
  \subref{fig:iter_pd_prog2ex_rev} the number of iterations of the 
  proposed method; 
  \subref{fig:iter_soc_prog2ex_rev} the number of systems of linear 
  equations solved in ReSNA; 
  \subref{fig:resid_eq_prog2ex_rev} 
  the residual in \eqref{eq:residual.1}; 
  \subref{fig:resid_compl_prog2ex_rev} 
  the residual in \eqref{eq:residual.2}; and 
  \subref{fig:resid_gap_ineq_prog2ex_rev} 
  the residual in \eqref{eq:residual.3}. 
  }
  \label{fig:prog2ex.result}
\end{figure*}

Consider the problem outlined in \reffig{fig:ex_body_2d}, where an 
elastic body is in the plane-stress state. 
The body is discretized uniformly as 
$N_{X} \times N_{Y}$ four-node quadrilateral (Q4) elements, 
where $N_{X} = 2.5N_{Y}$ and the value of $N_{Y}$ is varied to change 
the size of problem instance. 
As for implementation of the finite element method, 
we use a Matlab code due to \citet{ACSLS11}. 
Uniformly distributed vertical traction of $0.01$ 
is applied at the top edge. 
The bottom nodes can possibly contact with a flat rigid obstacle. 
Hence, the number of contact candidate nodes is $c = N_{X}$. 
The initial gaps are set to $g_{j}=0.01$ $(j=1,\dots,c)$. 
The coefficient of friction is $\mu=0.5$. 

\reffig{fig:prog2ex.result} collects the computational results for 
$N_{Y} = 26$, 40, 54, 68, 80, 92, 104, 114, 124, 134, and 144. 
Here, ``{\#}DOF'' means the number 
of degrees of freedom of the nodal displacements, i.e., $d$. 
\reffig{fig:time_log_dof_prog2ex_rev} shows the computation time. 
The dashed lines correspond to $O(d)$, $O(d^{2})$, and $O(d^{3})$. 
The proposed method clearly outperforms ReSNA in terms of computation time. 
For example, $382.6\,\mathrm{s}$ and $98222.0\,\mathrm{s}$ 
were required by the proposed method and ReSNA, respectively, to solve 
the instance with $d=54600$ and $c=260$. 
At the equilibrium state, about 40\% contact candidate nodes are free, 
10\% are in sliding contact, and 50\% are in sticking contact. 

\reffig{fig:iter_pd_prog2ex_rev} reports the number of iterations 
required by the proposed method. 
ReSNA solves a system of linear equations at step~2.1 of  
Algorithm~2 in \citep{HYF05}, to obtain a search direction. 
\reffig{fig:iter_soc_prog2ex_rev} reports the number of 
systems of linear equations solved at step~2.1. 
Increase of the number in \reffig{fig:iter_soc_prog2ex_rev} is moderate 
compared with \reffig{fig:iter_pd_prog2ex_rev}. 
Therefore, increase of computation time of ReSNA observed in 
\reffig{fig:time_log_dof_prog2ex_rev} may be due to increase of 
computation time for numerical solution of 
a system of linear equations  at each iteration. 

\reffig{fig:resid_eq_prog2ex_rev}, 
\reffig{fig:resid_compl_prog2ex_rev}, and 
\reffig{fig:resid_gap_ineq_prog2ex_rev} 
compare the residuals defined by 
\eqref{eq:residual.1}, \eqref{eq:residual.2}, and \eqref{eq:residual.3}, 
respectively.  
Concerning the residuals in \eqref{eq:residual.1} and 
\eqref{eq:residual.2}, it is observed in 
\reffig{fig:resid_eq_prog2ex_rev} and 
\reffig{fig:resid_compl_prog2ex_rev} that the solution obtained by the 
proposed method has smaller residuals for every problem instance. 
In contrast, in terms of the residual in \eqref{eq:residual.3}, the 
solution obtained by ReSNA has a smaller residual, although a residual 
of the solution obtained by the proposed method is also sufficiently 
small, as observed in \reffig{fig:resid_gap_ineq_prog2ex_rev}.

\subsection{Example (II): three-dimensional problem}
\label{sec:ex.3D}

\begin{figure*}[tbp]
  \centering
  \includegraphics[scale=0.45]{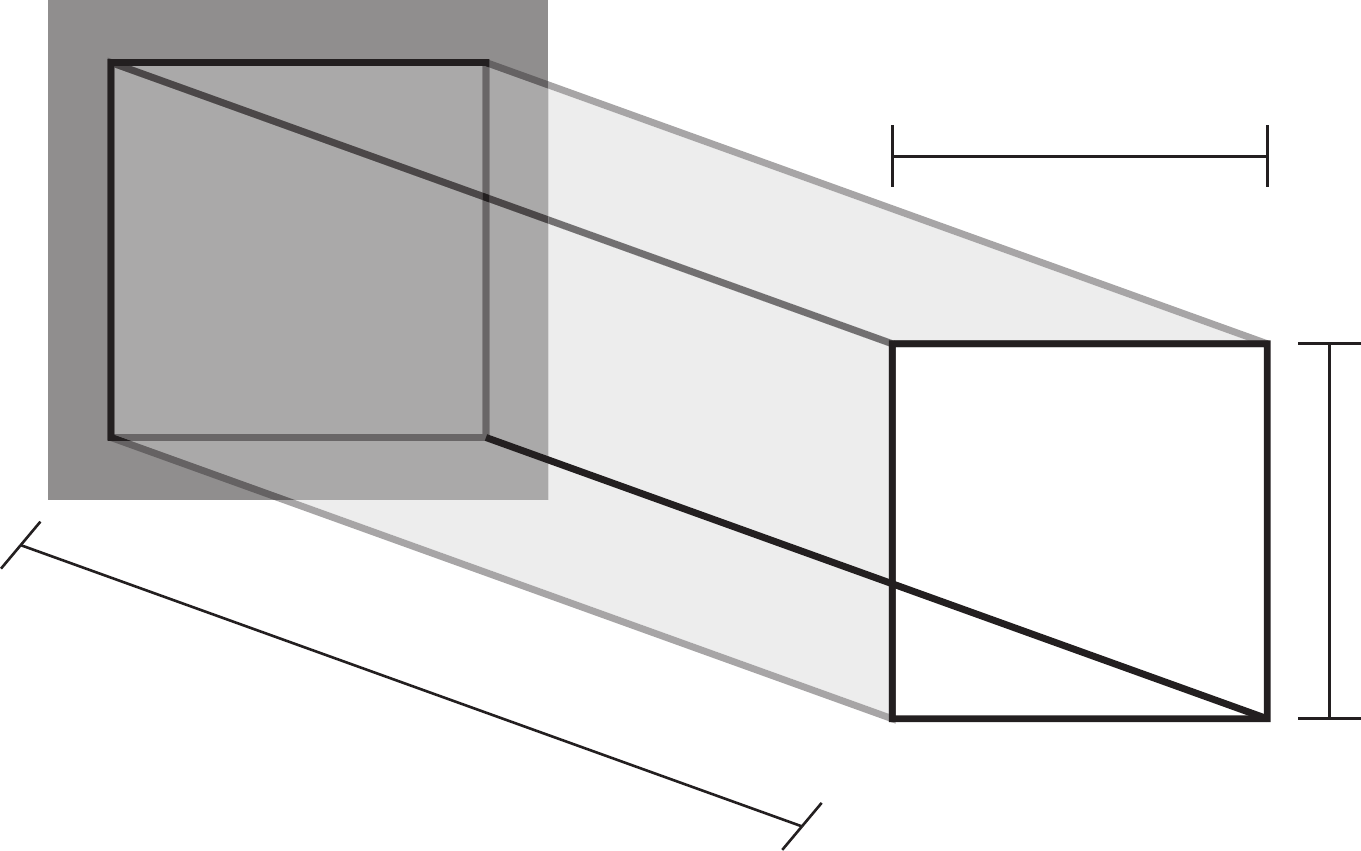}
  \begin{picture}(0,0)
    \put(-200,-150){
    \put(60,162){{\small $N_{X}$}}
    \put(194,190){{\small $N_{Z}$}}
    \put(154,244){{\small $N_{Y}$}}
    }
  \end{picture}
  \caption[]{Problem setting of example (II). }
  \label{fig:ex_body_3d}
\end{figure*}

\begin{figure*}[tbp]
  \centering
  \subfloat[]{
  \label{fig:time_log_dof_prog4ex_rev}
  \includegraphics[scale=0.50]{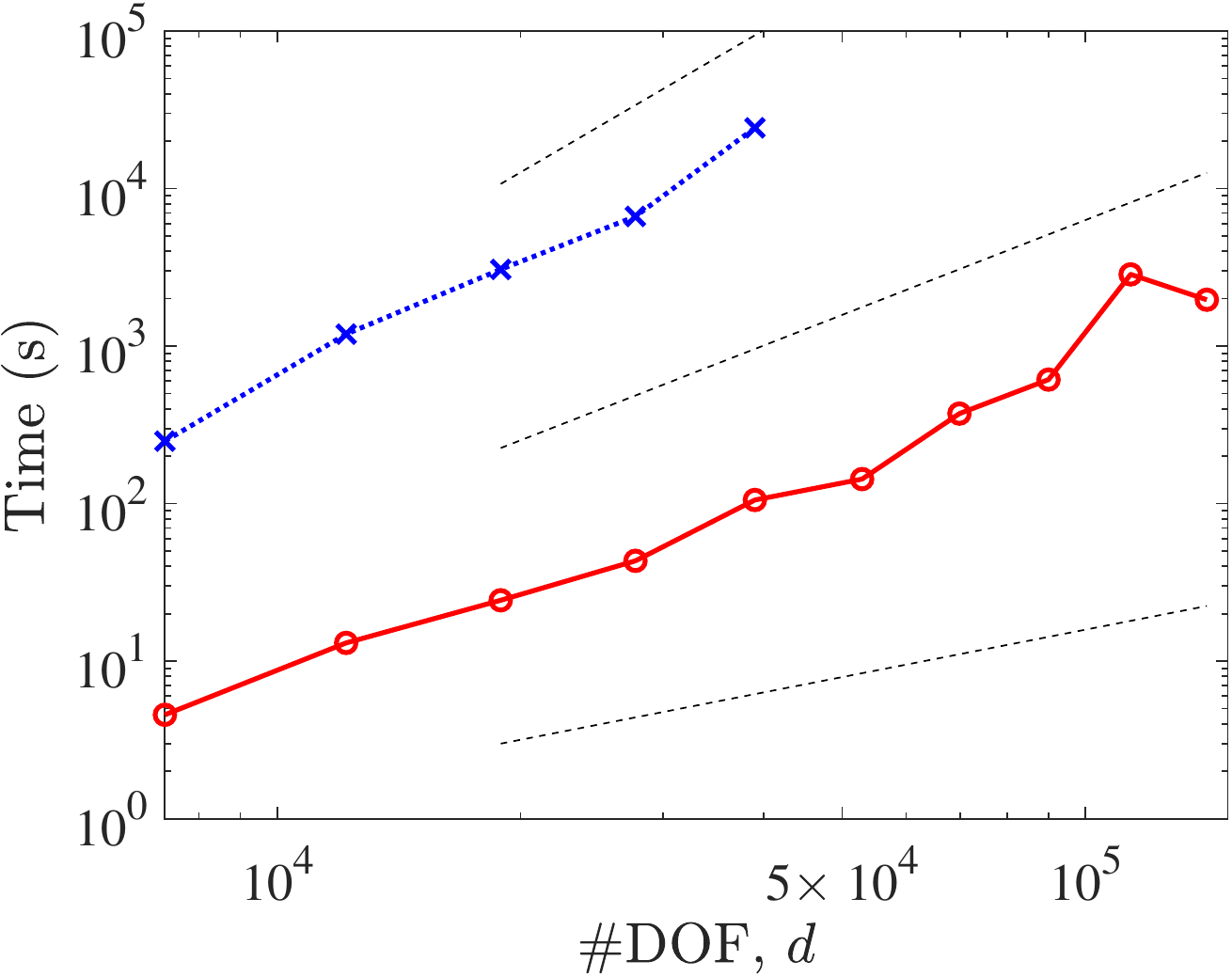}
  }
  \hfill
  \subfloat[]{
  \label{fig:iter_pd_prog4ex_rev}
  \includegraphics[scale=0.50]{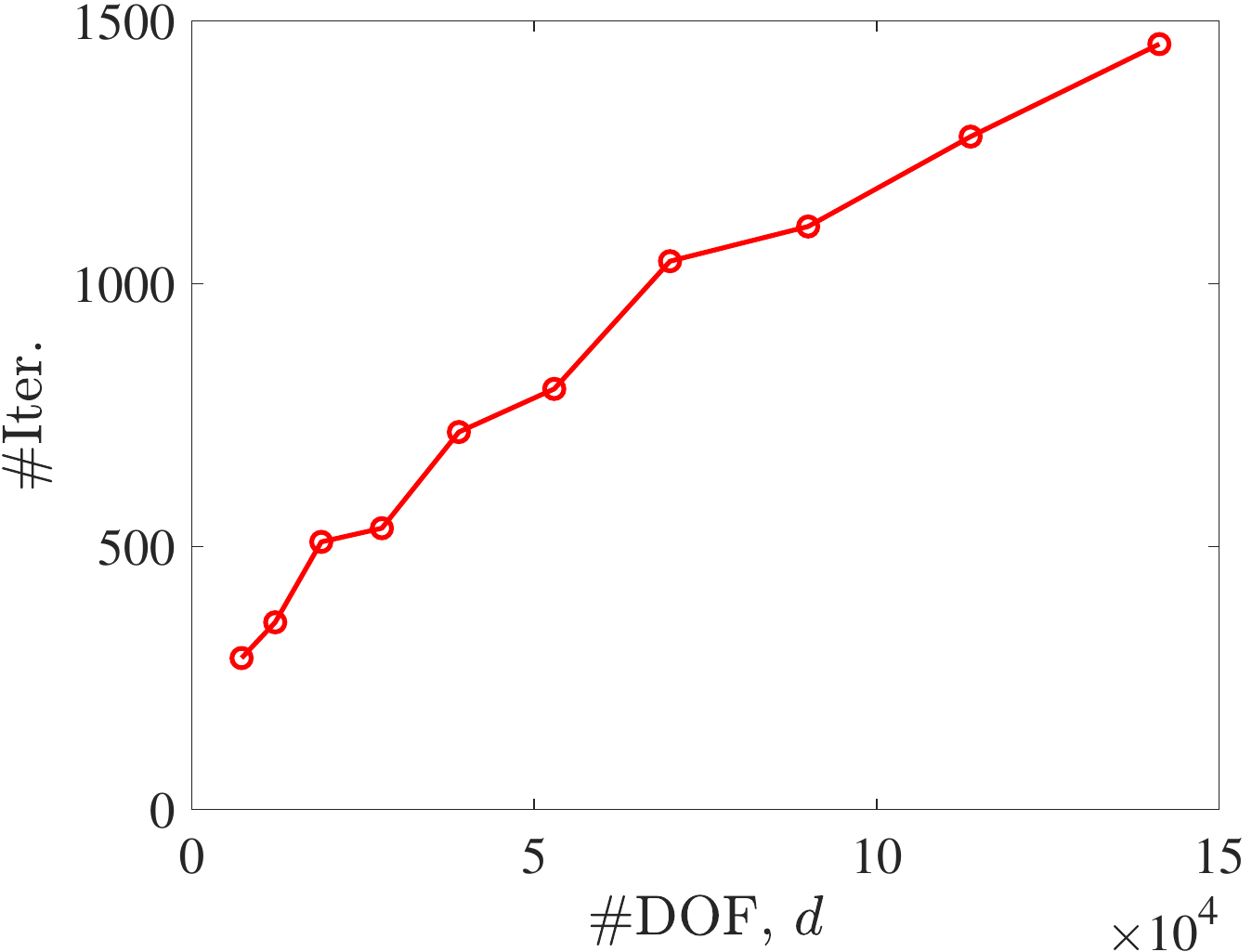}
  }
  \par
  \subfloat[]{
  \label{fig:iter_soc_prog4ex_rev}
  \includegraphics[scale=0.50]{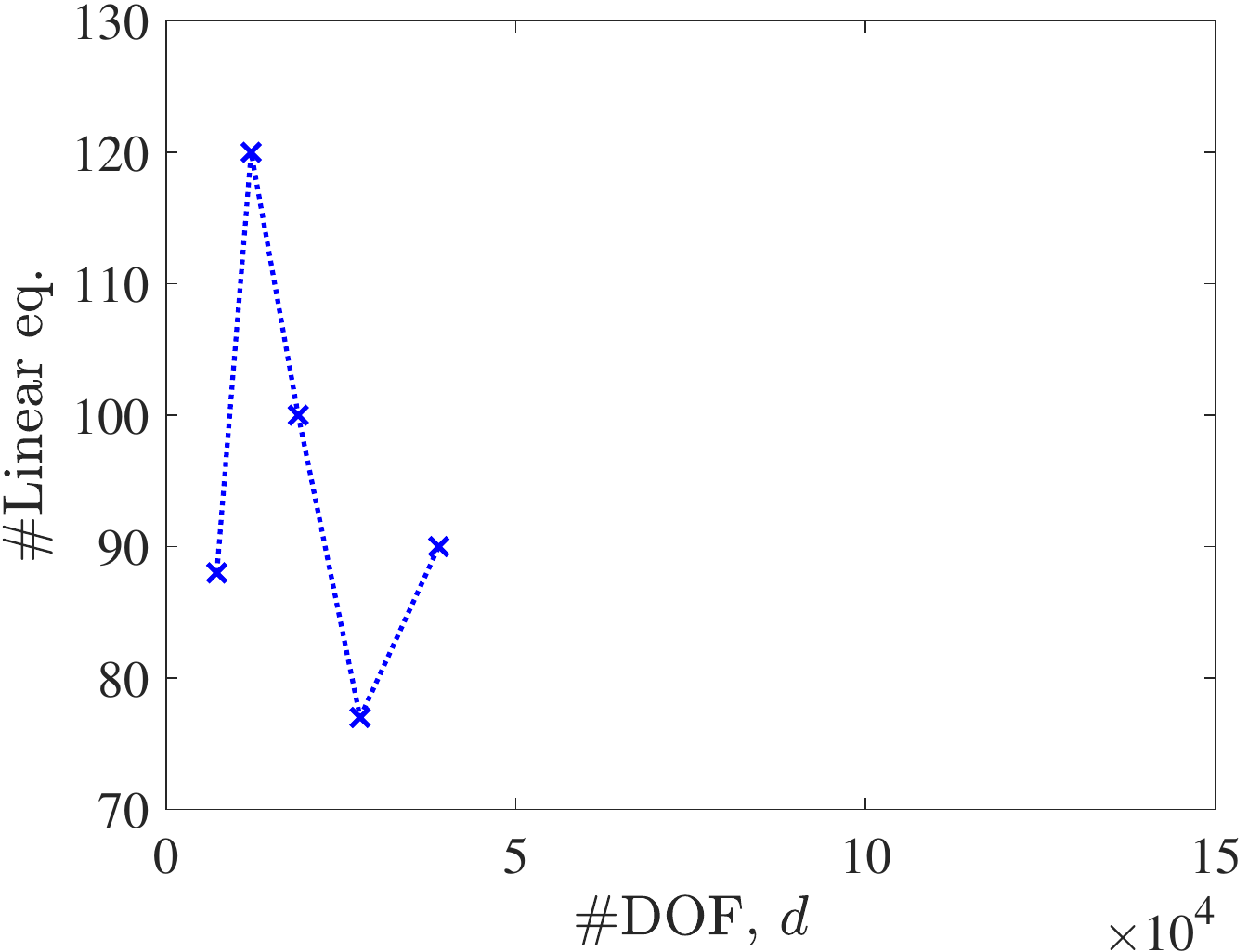}
  }
  \hfill
  \subfloat[]{
  \label{fig:resid_eq_prog4ex_rev}
  \includegraphics[scale=0.50]{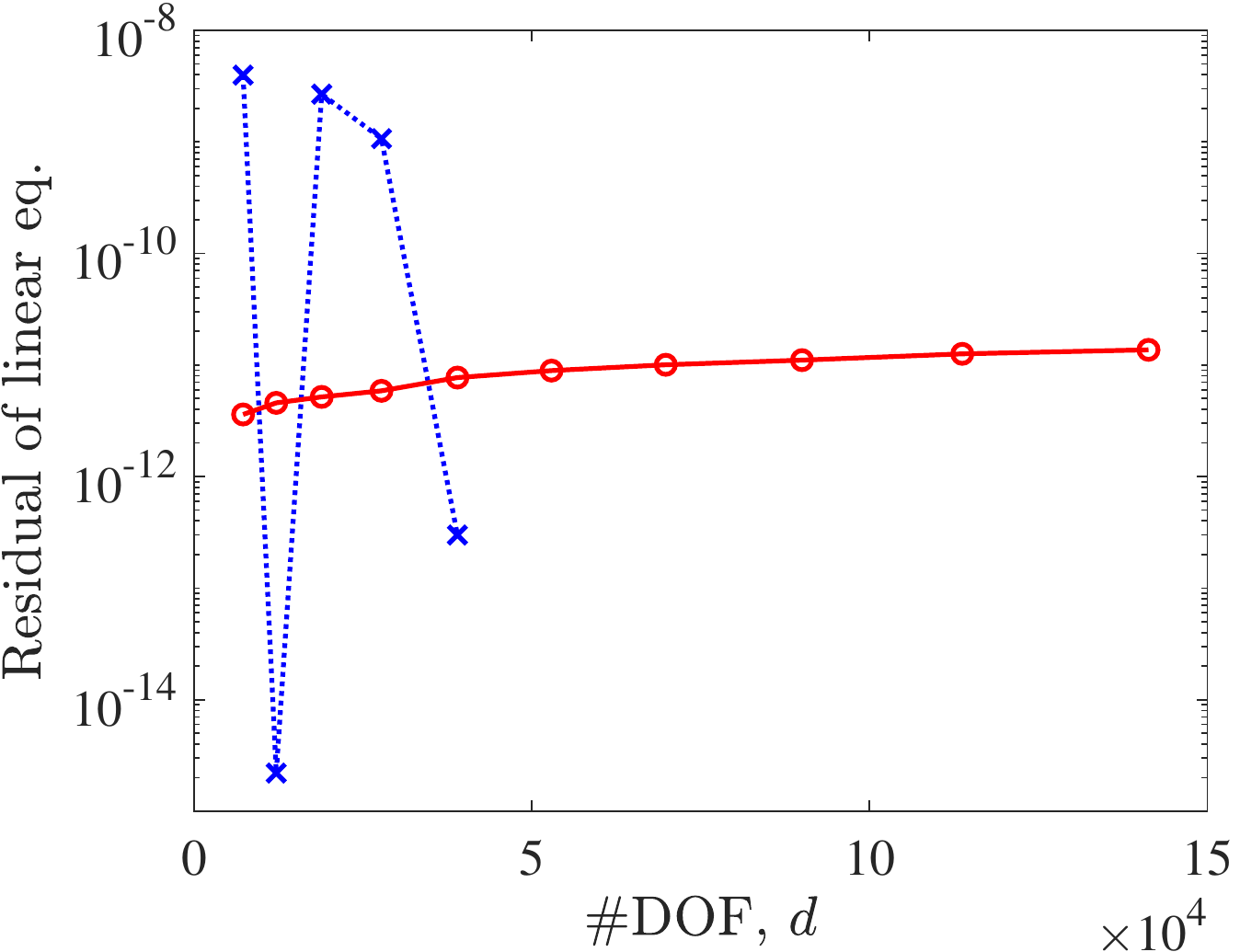}
  }
  \par
  \subfloat[]{
  \label{fig:resid_compl_prog4ex_rev}
  \includegraphics[scale=0.50]{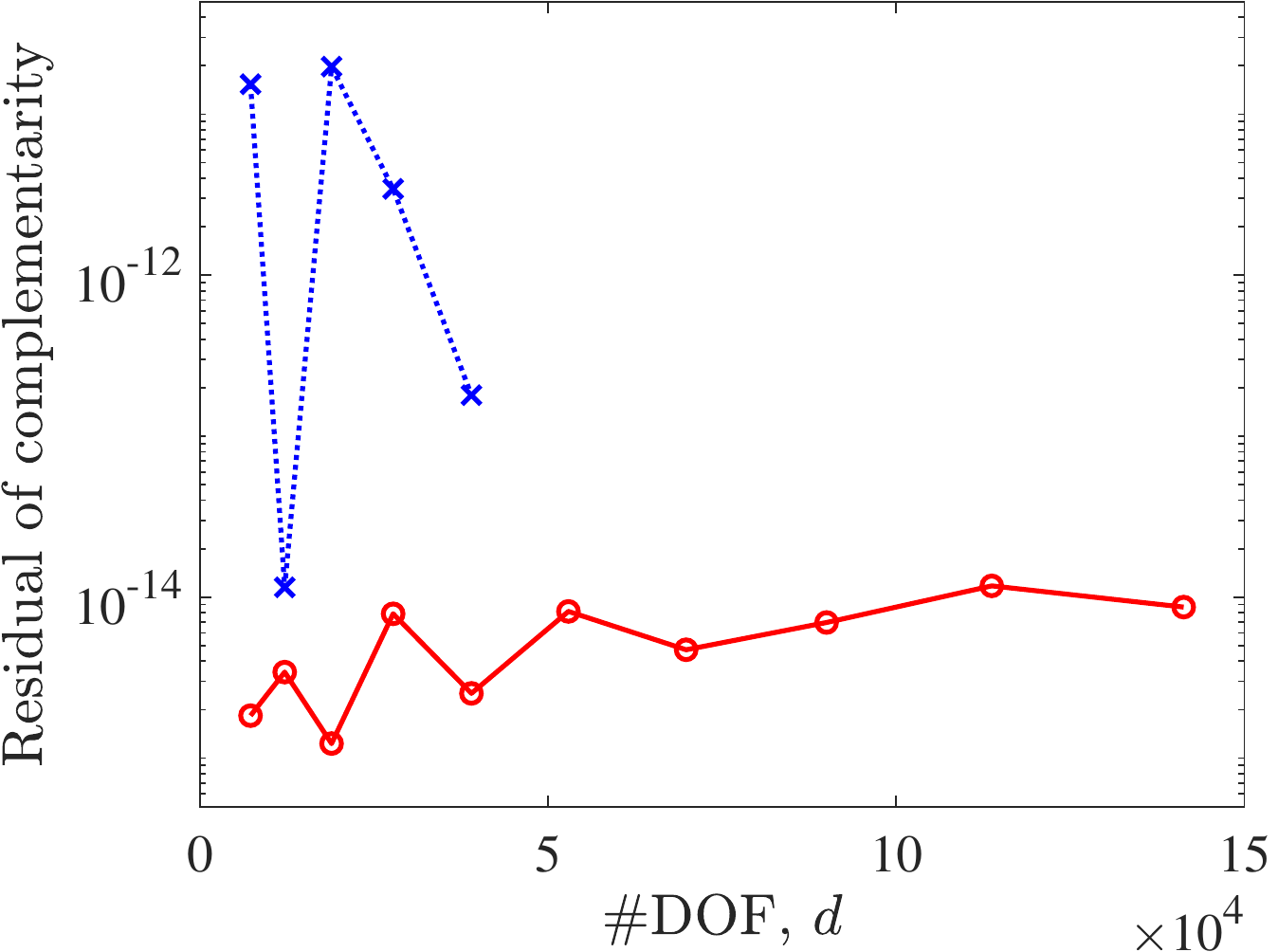}
  }
  \hfill
  \subfloat[]{
  \label{fig:resid_gap_ineq_prog4ex_rev}
  \includegraphics[scale=0.50]{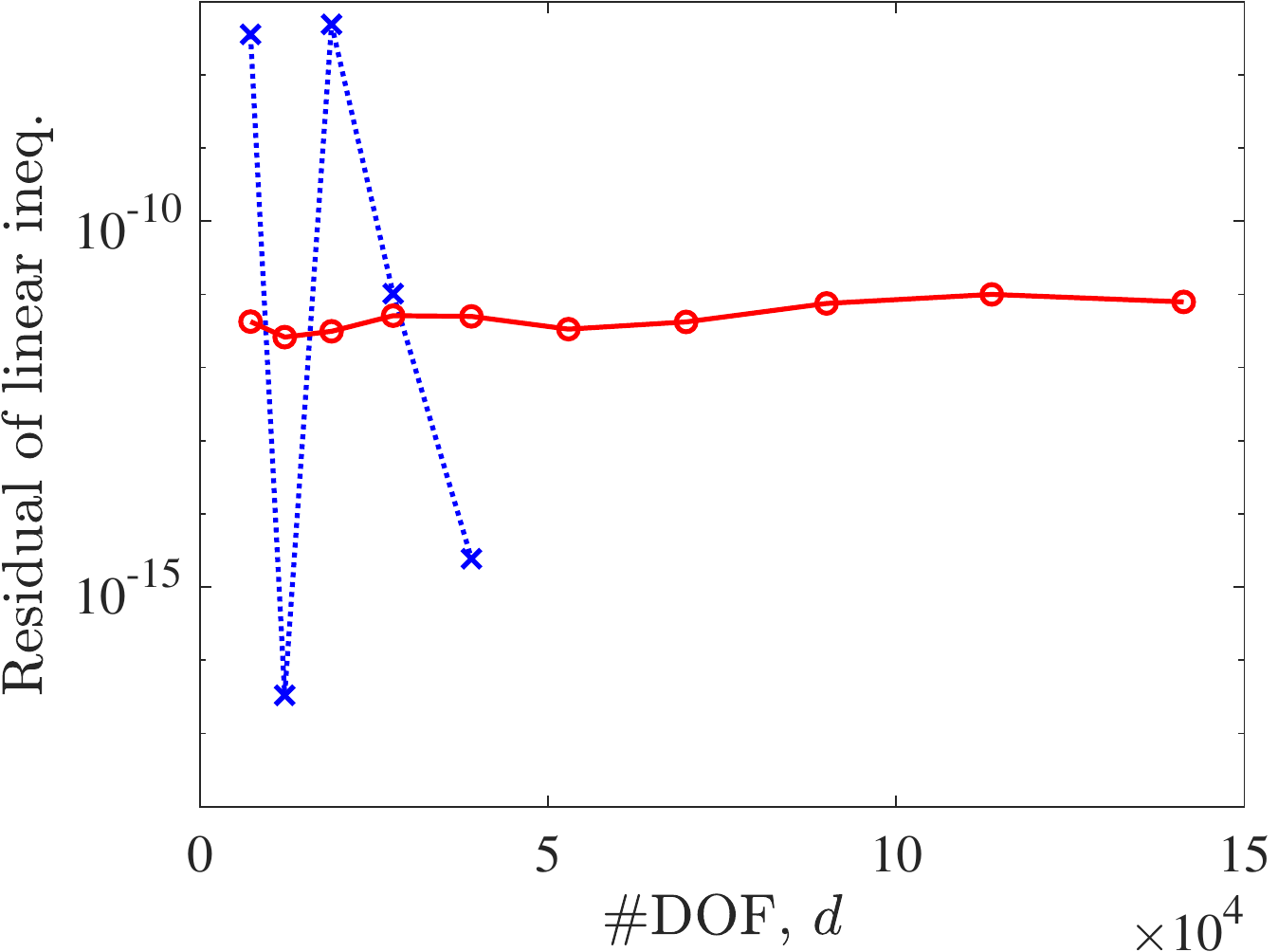}
  }
  \caption[]{Computational results of example (II) with $\mu=0.5$. 
  ``{\em Solid line\/}'' The proposed method; 
  and ``{\em dotted line\/}'' ReSNA. 
  \subref{fig:time_log_dof_prog4ex_rev} The computation time; 
  \subref{fig:iter_pd_prog4ex_rev} the number of iterations of the 
  proposed method; 
  \subref{fig:iter_soc_prog4ex_rev} the number of systems of linear 
  equations solved in ReSNA; 
  \subref{fig:resid_eq_prog4ex_rev} 
  the residual in \eqref{eq:residual.1}; 
  \subref{fig:resid_compl_prog4ex_rev} 
  the residual in \eqref{eq:residual.2}; and 
  \subref{fig:resid_gap_ineq_prog4ex_rev} 
  the residual in \eqref{eq:residual.3}. 
  }
  \label{fig:prog_4.mu0p5}
\end{figure*}

\begin{figure*}[tbp]
  \centering
  \subfloat[]{
  \label{fig:prog4_pcg_iter}
  \includegraphics[scale=0.50]{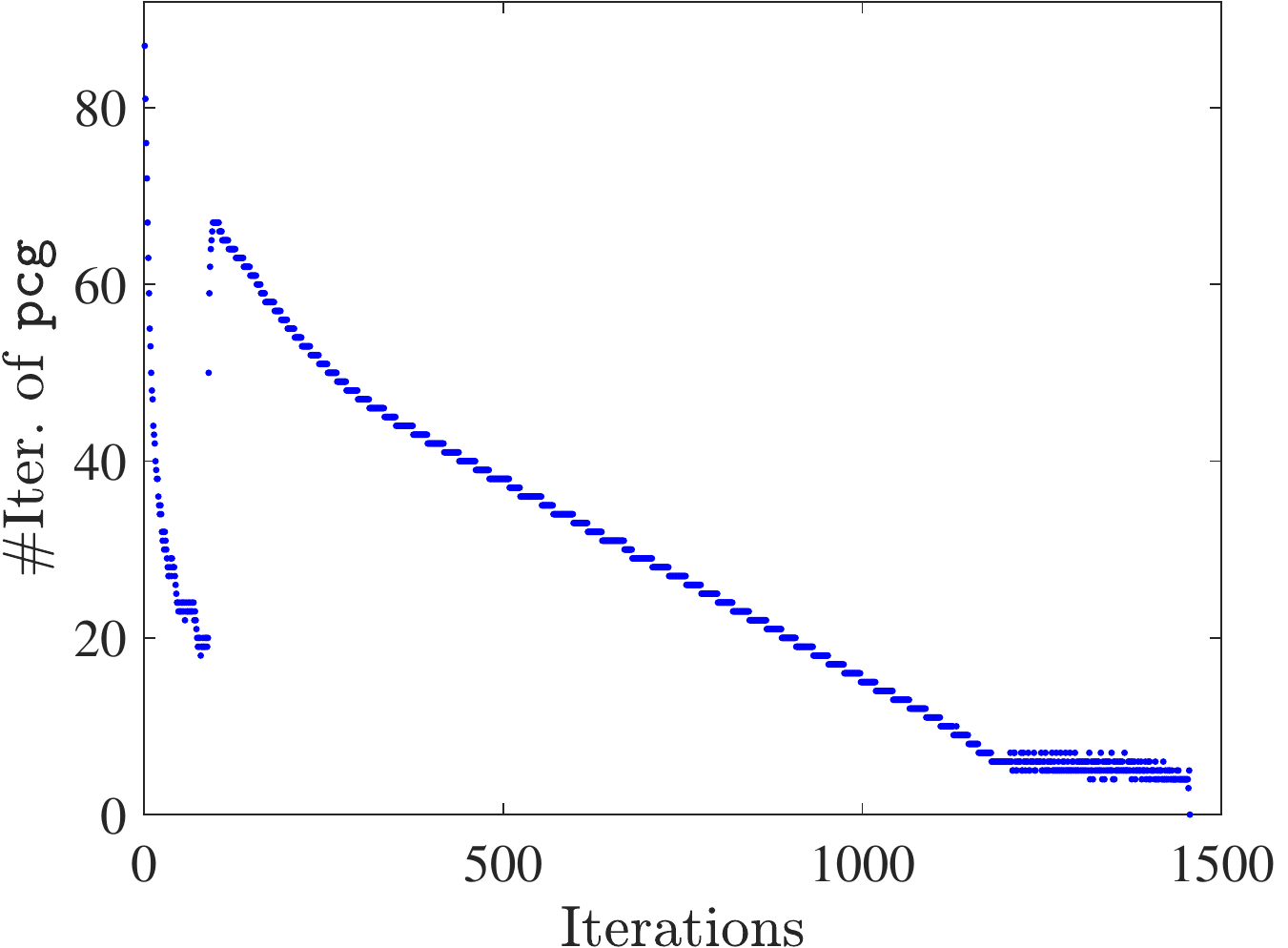}
  }
  \hfill
  \subfloat[]{
  \label{fig:prog4_resid_force_eq}
  \includegraphics[scale=0.50]{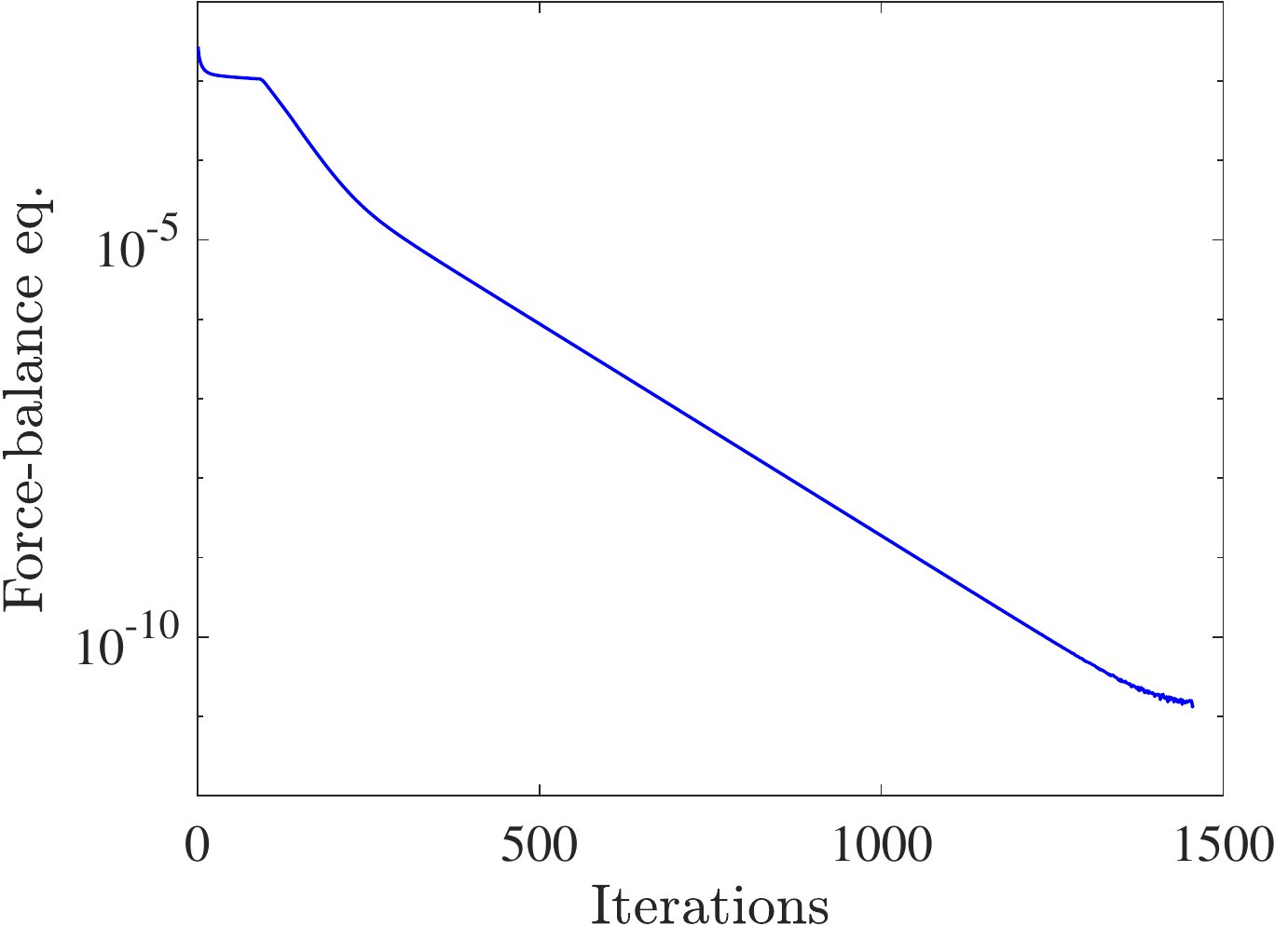}
  }
  \par
  \subfloat[]{
  \label{fig:prog4_resid_comp}
  \includegraphics[scale=0.50]{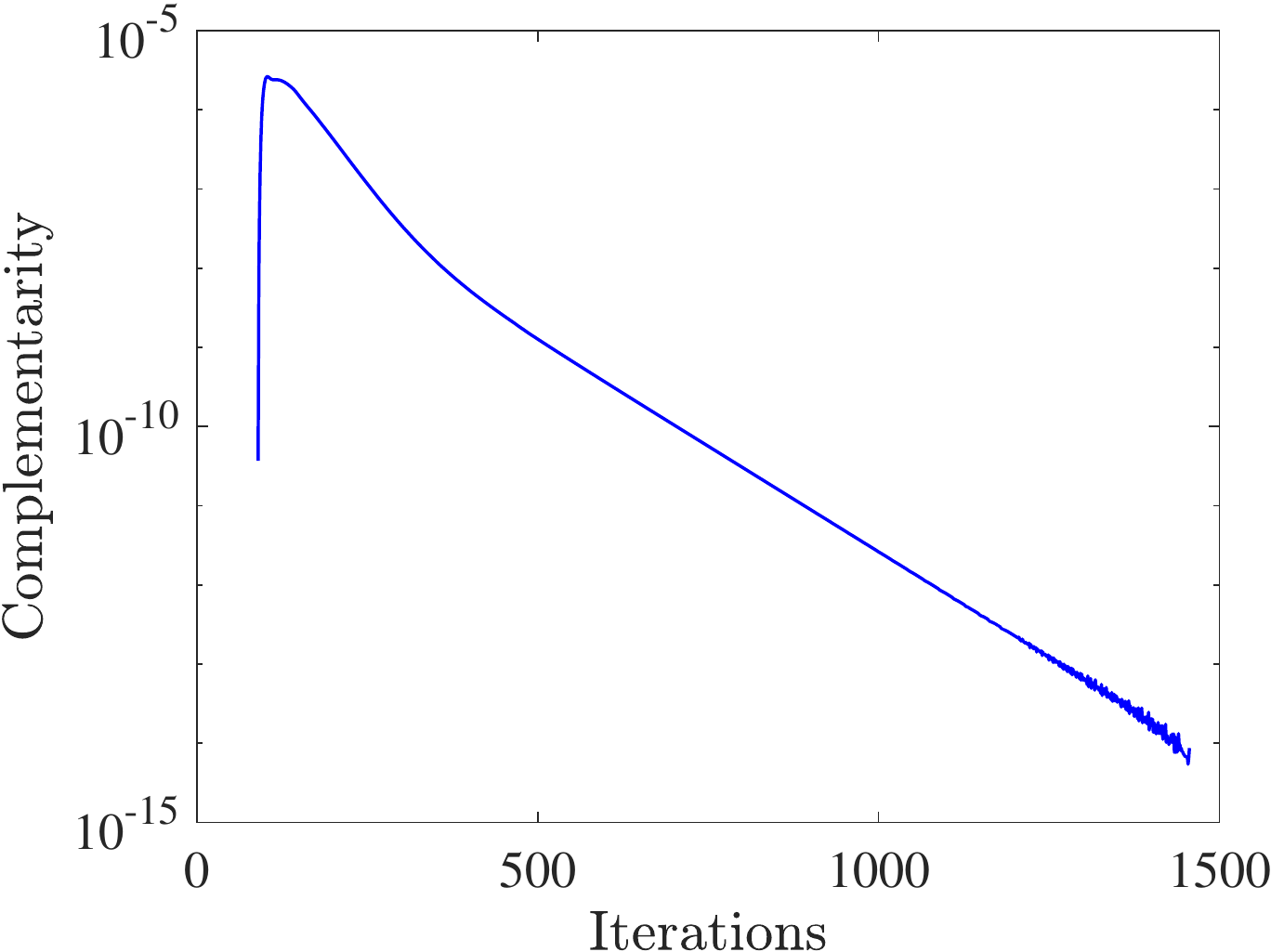}
  }
  \hfill
  \subfloat[]{
  \label{fig:prog4_resid_gap_ineq}
  \includegraphics[scale=0.50]{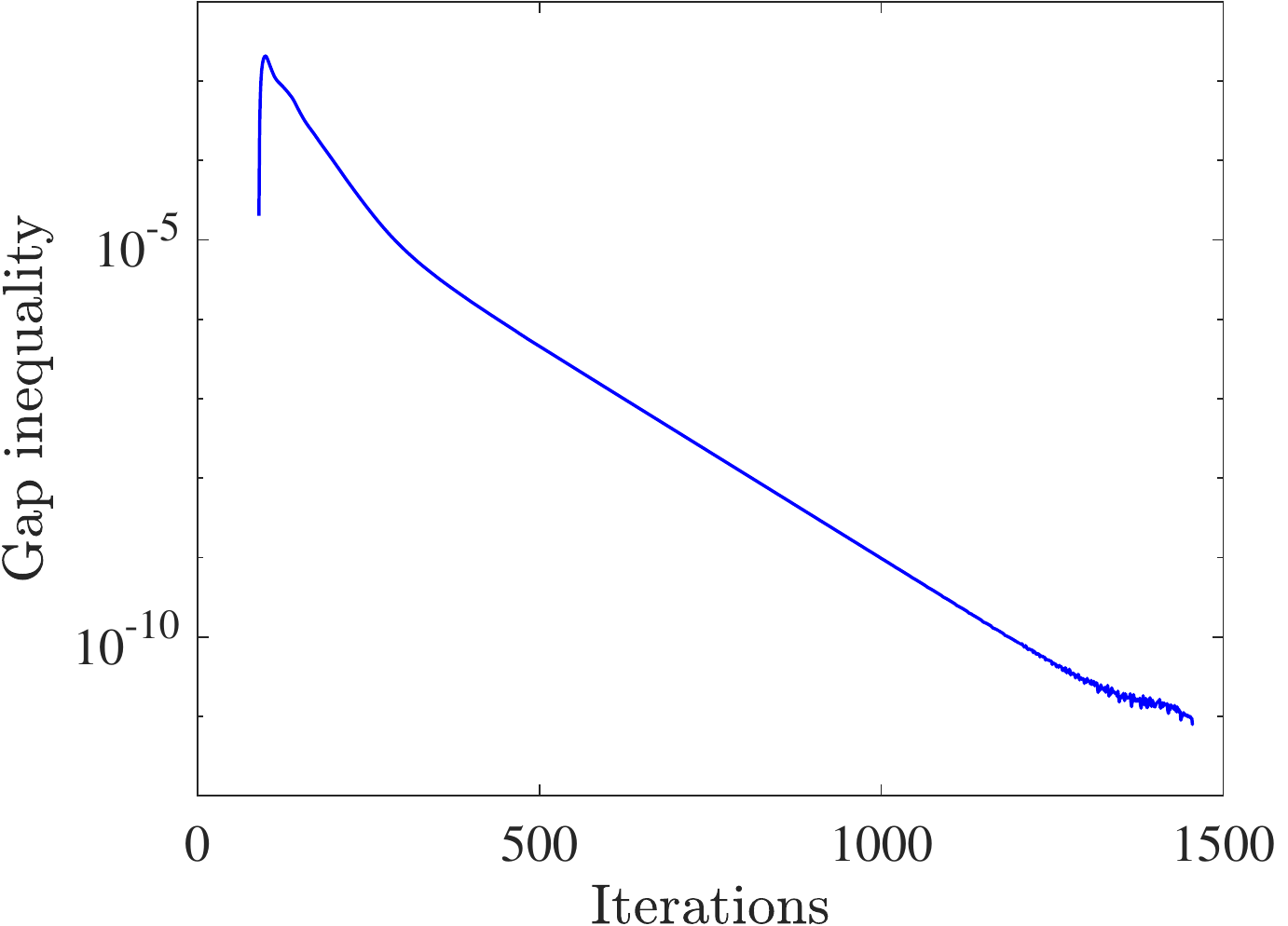}
  }
  \caption[]{Iteration history of the proposed method for example (II) 
  with $\mu=0.5$, $N_{X}=56$, $N_{Y}=N_{Z}=28$, $d=141288$, and $c=1624$. 
  \subref{fig:prog4_pcg_iter} The number of iterations of the conjugate 
  gradient method; 
  \subref{fig:prog4_resid_force_eq} 
  the residual in \eqref{eq:residual.1}; 
  \subref{fig:prog4_resid_comp} 
  the residual in \eqref{eq:residual.2}; and 
  \subref{fig:prog4_resid_gap_ineq} 
  the residual in \eqref{eq:residual.3}. 
  }
  \label{fig:prog4.iteration}
\end{figure*}

\begin{figure*}[tbp]
  \centering
  \subfloat[]{
  \label{fig:time_log_dof_prog4ex_fric}
  \includegraphics[scale=0.50]{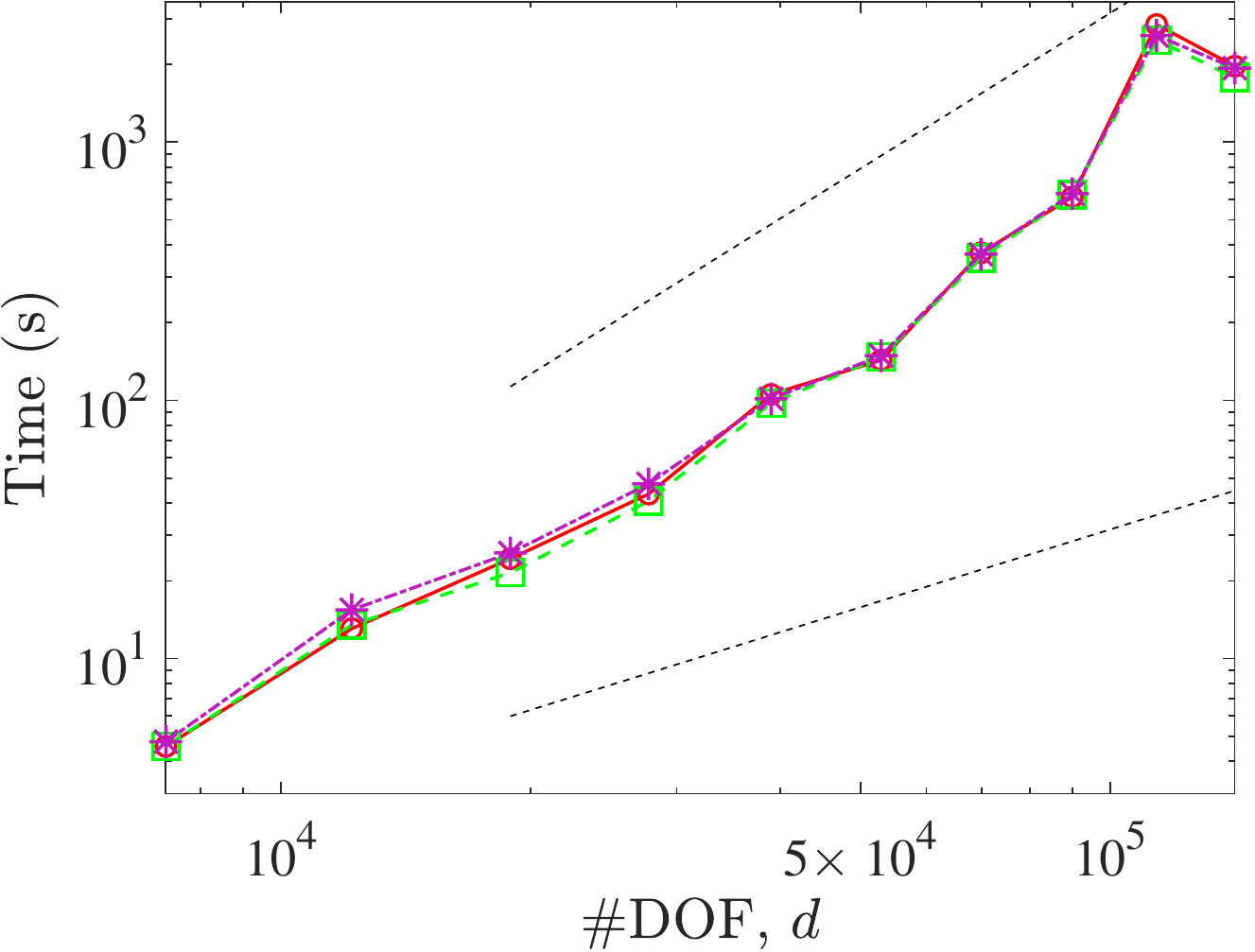}
  }
  \hfill
  \subfloat[]{
  \label{fig:resid_eq_prog4ex_fric}
  \includegraphics[scale=0.50]{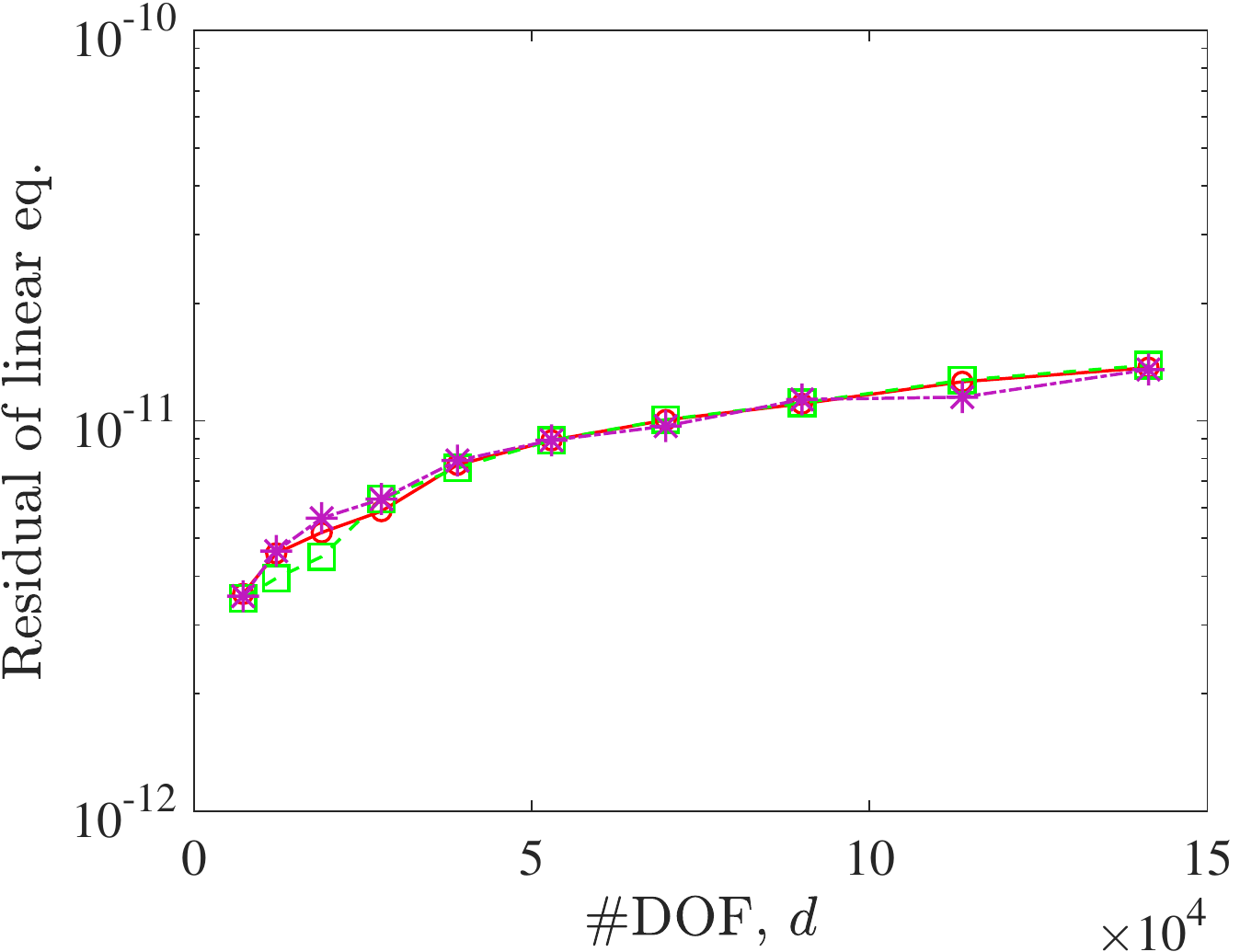}
  }
  \par
  \subfloat[]{
  \label{fig:resid_compl_prog4ex_fric}
  \includegraphics[scale=0.50]{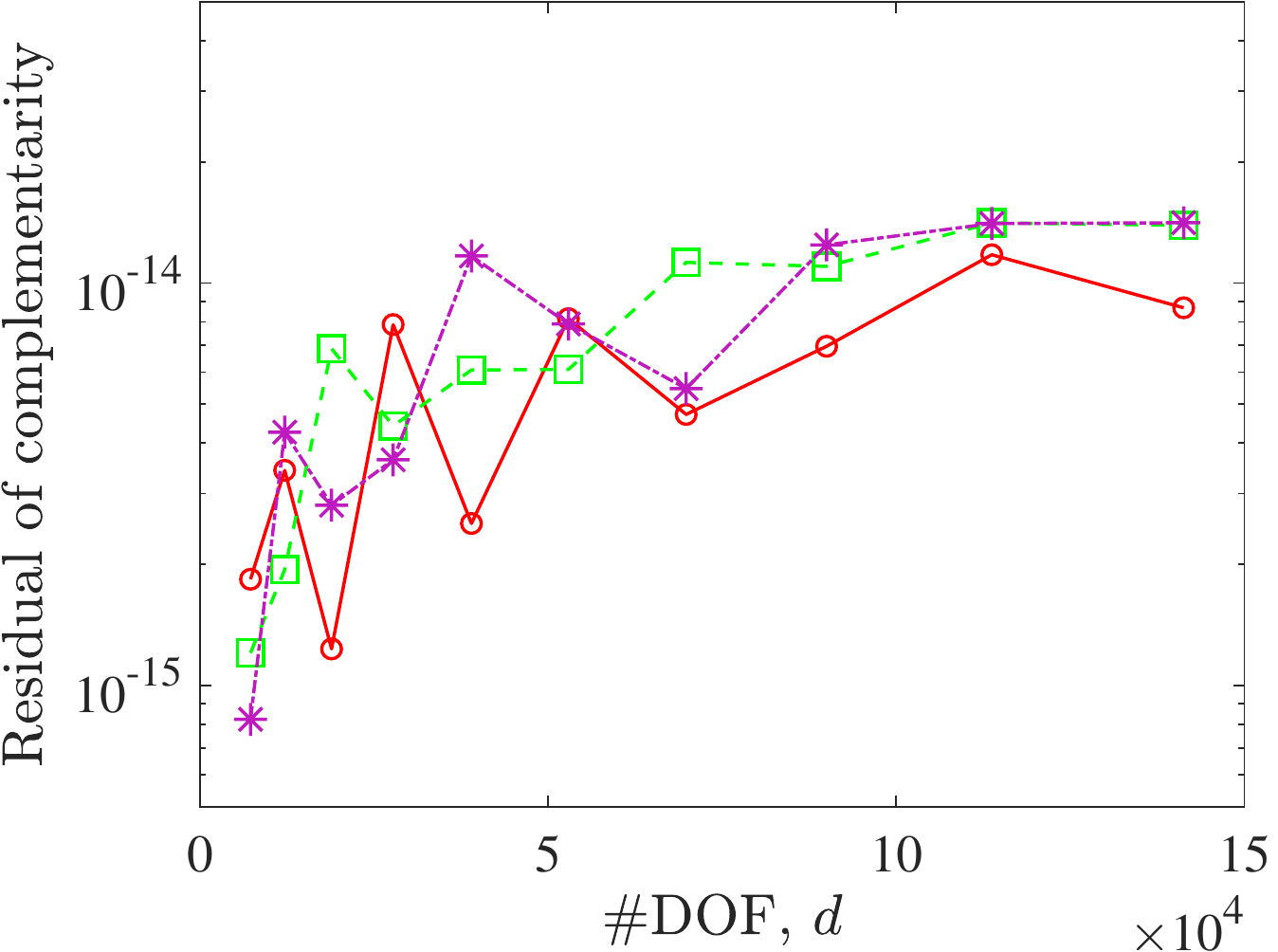}
  }
  \hfill
  \subfloat[]{
  \label{fig:resid_gap_ineq_prog4ex_fric}
  \includegraphics[scale=0.50]{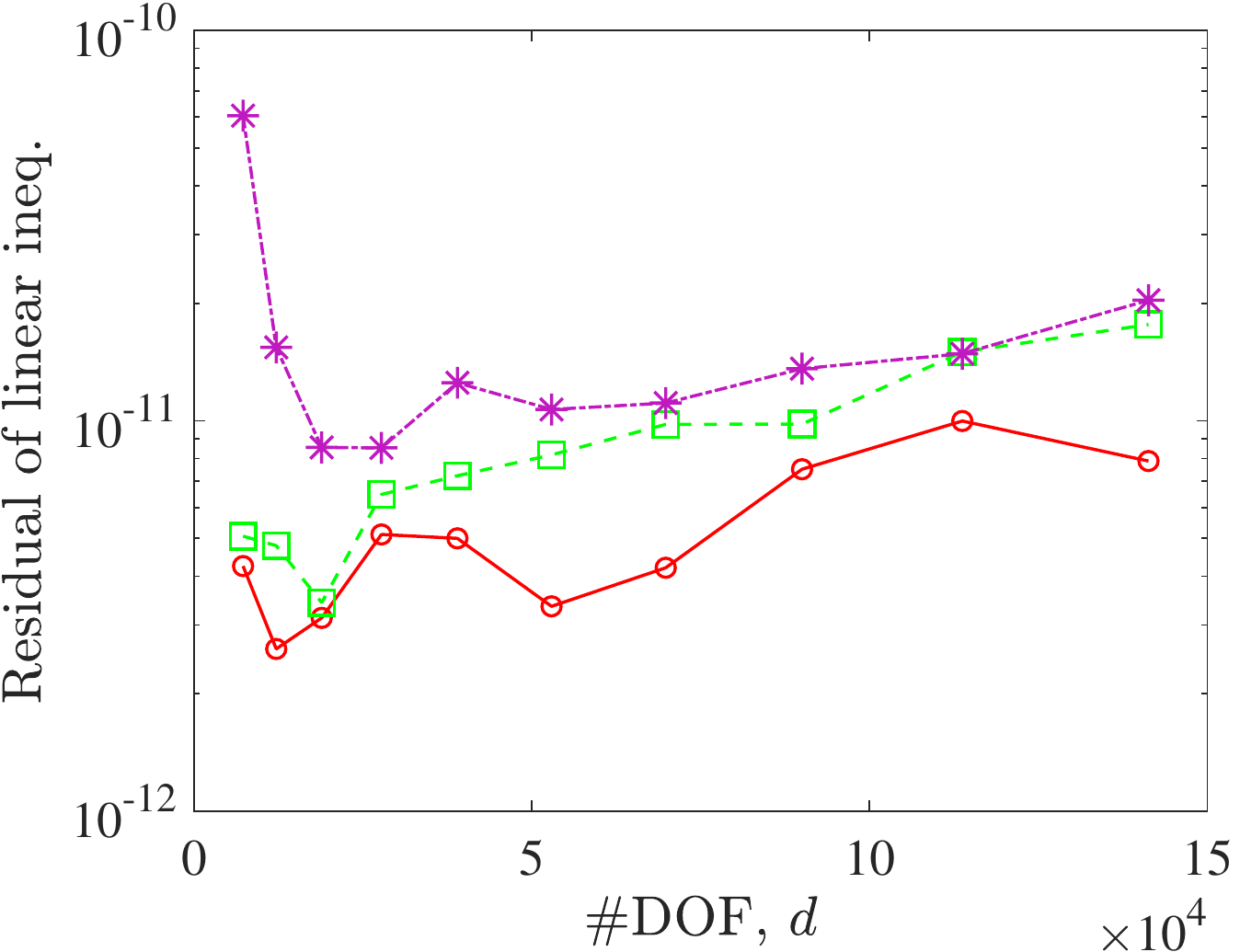}
  }
  \caption[]{Computational results of the proposed method for example 
  (II) with three different values of the friction coefficient. 
  ``{\em Solid line\/}'' $\mu=0.5$; 
  ``{\em dashed line\/}'' $\mu=1.0$; and 
  ``{\em dashed-doted line\/}'' $\mu=1.5$. 
  \subref{fig:time_log_dof_prog4ex_fric} The computational time; 
  \subref{fig:resid_eq_prog4ex_fric} 
  the residual in \eqref{eq:residual.1}; 
  \subref{fig:resid_compl_prog4ex_fric} 
  the residual in \eqref{eq:residual.2}; and 
  \subref{fig:resid_gap_ineq_prog4ex_fric} 
  the residual in \eqref{eq:residual.3}. 
  }
  \label{fig:prog4ex_fric}
\end{figure*}

Consider a three-dimensional version of example (I) outlined in 
\reffig{fig:ex_body_3d}. 
The elastic body is discretized uniformly as 
$N_{X} \times N_{Y} \times N_{Z}$ 8-node hexahedron elements, 
where $N_{X} = 2N_{Y} = 2N_{Z}$. 
A Matlab code due to \citet{FS20} is used as for 
implementation of the finite element method. 
Uniformly distributed vertical traction of $5 \times 10^{-3}$ 
is applied at all the top nodes. 
The bottom nodes are considered as contact candidate nodes, 
the number of which is $c = N_{X}(N_{Y}+1)$. 
The initial gaps are $g_{j}=0.005$ $(j=1,\dots,c)$.
At the equilibrium state, about 33\% contact candidate nodes are free, 
11\% are in sliding contact, and 56\% are in sticking contact. 

\reffig{fig:prog_4.mu0p5} reports the computational results with 
$\mu=0.5$, where $N_{Y} = 10$, 12, 14, 16, 18, 20, 22, 24, 26, and 28. 
It is observed in \reffig{fig:time_log_dof_prog4ex_rev}  that the proposed 
method outperforms ReSNA from the viewpoint of computation time. 
Concerning the number of iterations, 
\reffig{fig:iter_pd_prog4ex_rev} and \reffig{fig:iter_soc_prog4ex_rev} 
show trends similar to the ones observed in example (I). 
From \reffig{fig:resid_eq_prog4ex_rev}, 
\reffig{fig:resid_compl_prog4ex_rev}, and 
\reffig{fig:resid_gap_ineq_prog4ex_rev} we see that 
accuracy of the solutions obtained by the proposed 
method is comparable with that of the solutions obtained by ReSNA. 

\reffig{fig:prog4.iteration} shows the iteration history of the proposed 
method. 
The number of iterations required by \texttt{pcg} for solving a system of 
linear equations in line~\ref{alg:frictional.accelerated.equation} 
of \refalg{alg:frictional.accelerated.2} decreases as the algorithm 
approaches to termination. 
The residual in \eqref{eq:residual.1} decreases monotonically. 
In contrast, the residuals in \eqref{eq:residual.2} and 
\eqref{eq:residual.3} are equal to zero for about the first 90 iterations. 
This is because the contact candidates nodes are free and have no 
reactions at these iterations as we set 
$\Delta\bi{u}^{(0)}=\bi{0}$ and $\bi{r}^{(0)} = \bi{0}$. 
The residuals become positive when the displacement violates the 
non-penetration condition, and then decrease monotonically. 

For $\mu=0.5$, $1$, and $1.5$, \reffig{fig:prog4ex_fric} compares 
performance of the proposed method. 
We can observe that the performance is irrelevant to a value of 
the friction coefficient.

\section{Conclusions}
\label{sec:conclusion}

This paper has developed a fast fist-order optimization-based method for 
a quasi-static contact problem with Coulomb's friction. 
The method is designed based on an accelerated primal-dual algorithm 
solving a convex optimization problem that approximates the contact problem. 

In the numerical experiments on problem instances with up to 140 
thousands degrees of freedom of the nodal displacements, it has been 
observed that the proposed method successfully converges to a solution 
for every problem instance, although the method has no guarantee of 
convergence. 
Also, it has been demonstrated that the proposed method outperforms 
a regularized and smoothed Newton method for the second-order cone 
complementarity problem. 
Furthermore, the proposed method is easy to implement.

\paragraph{Acknowledgments}

This work is supported by 
JSPS KAKENHI 17K06633, 21K04351, and 
JST CREST Grant No.~JPMJCR1911, Japan.

\appendix

\section{Formulation as second-order cone linear complementarity problem}
\label{sec:linear_form}

ReSNA \citep{Hay20} is a Matlab software package implementing a 
regularized smoothing Newton method proposed by \citet{HYF05}. 
It solves a second-order cone linear complementarity problem (SOCLCP) 
in the following form: 
\begin{subequations}\label{P.linear.SOC.1}%
  \begin{align}
    & \KC \ni \bi{x} \perp \bi{y} \in \KC , 
    \label{P.linear.SOC.1.1} \\
    & \bi{y} = M_{11} \bi{x} + M_{12} \bi{v} + \bi{w}_{1} , \\
    & M_{21} \bi{x} + M_{22} \bi{v} + \bi{w}_{2} = \bi{0} . 
  \end{align}
\end{subequations}
Here, $\bi{x}$, $\bi{y}$, and $\bi{v}$ are unknown variables, 
and $\KC$ is a Cartesian product of some second-order cones. 
In the numerical experiments reported in section~\ref{sec:ex}, we use 
ReSNA for comparison. 

The method proposed in this paper attempts to solve problem 
\eqref{eq:nonlinear.complementarity} in section~\ref{sec:problem}.  
\citet{KMPc06} showed that this problem can be recast as the 
following SOCLCP: 
\begin{subequations}\label{P.linear.SOC.2}%
  \begin{alignat}{3}
    & K \Delta\bi{u} =
    \bi{p} + T_{\rr{n}} \bi{r}_{\rr{n}} + T_{\rr{t}} \bi{r}_{\rr{t}} , 
    \label{P.linear.SOC.2.1} \\
    & \Re_{+} \ni 
    (g_{j} - \bi{t}_{\rr{n}j}^{\top} \Delta\bi{u}) 
    \perp
    (-r_{\rr{n}j}) \in \Re_{+} , 
    \quad j=1,\dots,c , 
    \label{P.linear.SOC.2.2} \\
    & L^{3} \ni 
    \begin{bmatrix}
      \lambda_{j} \\
      T_{\rr{t}j}^{\top} \Delta\bi{u} \\
    \end{bmatrix}
    \perp 
    \begin{bmatrix}
      -\mu r_{\rr{n}j} \\
      \bi{r}_{\rr{t}j} \\
    \end{bmatrix}
    \in L^{3} , 
    \quad j=1,\dots,c . 
    \label{P.linear.SOC.2.3}
  \end{alignat}
\end{subequations}
Here, $\Re_{+}^{n}$ and $L^{n}$ are the nonnegative orthant 
and the second-order cone, respectively, i.e., 
\begin{align*}
  \Re_{+}^{n} &= \{
  \bi{x} \in \Re^{n} 
  \mid
  \bi{x} \ge \bi{0} \} , \\
  L^{n} &= \{
  (x_{0},\bi{x}_{1}) \in \Re \times \Re^{n-1} 
  \mid
  x_{0} \ge \| \bi{x}_{1} \| 
  \} , 
\end{align*}
and $\lambda_{j} \in \Re$ $(j=1,\dots,c)$ are additional unknown variables. 

Problem \eqref{P.linear.SOC.2} can be embedded into the form of 
\eqref{P.linear.SOC.1} as follows. 
Put 
\begin{align}
  \bi{x} = 
  \begin{bmatrix}
    \bi{g} - T_{\rr{n}}^{\top} \Delta\bi{u} \\
    \lambda_{1} \\
    T_{\rr{t}1}^{\top} \Delta\bi{u} \\
    \vdots \\
    \lambda_{c} \\
    T_{\rr{t}c}^{\top} \Delta\bi{u} \\
  \end{bmatrix}
  , \quad
  \bi{y} = 
  \begin{bmatrix}
    -\bi{r}_{\rr{n}} \\
    -\mu r_{\rr{n}1} \\
    \bi{r}_{\rr{t}1} \\
    \vdots \\
    -\mu r_{\rr{n}c} \\
    \bi{r}_{\rr{t}c} \\
  \end{bmatrix}
  , \quad
  \bi{v} = 
  \begin{bmatrix}
    \Delta\bi{u} \\
    \bi{\lambda} \\
    \bi{r}_{\rr{n}} \\
    \bi{r}_{\rr{t}}
  \end{bmatrix}
  \label{eq:def.x.y.v}
\end{align}
to see that \eqref{P.linear.SOC.2.2} and \eqref{P.linear.SOC.2.3} are 
equivalently rewritten as \eqref{P.linear.SOC.1.1} with 
$\KC \subset \Re^{4c}$ defined by 
\begin{align*}
  \KC 
  = \Re_{+}^{c} \times \LC ^{3} \times \dots \times \LC^{3} . 
\end{align*}
Define $\bi{e}_{1} \in \Re^{3}$ and $E_{2} \in \Re^{3 \times 2}$ by 
\begin{align*}
  \bi{e}_{1} = 
  \begin{bmatrix}
    1 \\ 0 \\ 0 \\
  \end{bmatrix}
  , \quad
  E_{2} = 
  \begin{bmatrix}
    0 & 0 \\ 1 & 0 \\ 0 & 1 \\
  \end{bmatrix}
  . 
\end{align*}
We see that the relation between $\bi{y}$ and $\bi{v}$ 
in \eqref{eq:def.x.y.v} can be written as 
\begin{align*}
  \bi{y} &= M_{12} \bi{v} 
\end{align*}
with $M_{12} \in \Re^{(4c) \times (d+4c)}$ defined by 
\begin{align}
  M_{12} = 
  \begin{bmatrix}
    O_{c,d} & O_{c,c} & -I_{c} & O_{c,2c} \\
    O_{3c,d} & O_{3c,c} & -\mu I_{c} \otimes \bi{e}_{1} 
    & I_{c} \otimes E_{2} \\
  \end{bmatrix}
  , 
  \label{eq:matrix.M.1.2}
\end{align}
where we use $\otimes$ to denote the Kronecker product. 
Similarly, define $M_{22}^{(1)} \in \Re^{(4c) \times (d+4c)}$ and 
$\bi{w}_{2}^{(1)} \in \Re^{4c}$ by 
\begin{align*}
  M_{22}^{(1)} = 
  \begin{bmatrix}
    T_{\rr{n}}^{\top} & O_{c,c} & O_{c,c} & O_{c,2c} \\
    -(I_{c} \otimes E_{2})T_{\rr{t}}^{\top} & -I_{c} \otimes \bi{e}_{1} 
    & O_{3c,c} & O_{3c,2c} \\
  \end{bmatrix}
  ,   \quad
  \bi{w}_{2}^{(1)} = 
  \begin{bmatrix}
    -\bi{g} \\ O_{3c,1} \\
  \end{bmatrix}
\end{align*}
to see that the relation between $\bi{x}$ and $\bi{v}$ is reduced to 
\begin{align*}
  \bi{x} + M_{22}^{(1)} \bi{v} + \bi{w}_{2}^{(1)}
  = \bi{0} . 
\end{align*}
Furthermore, define $M_{22}^{(2)} \in \Re^{d \times (d+4c)}$ and 
$\bi{w}_{2}^{(2)} \in \Re^{d}$ by 
\begin{align*}
  M_{22}^{(2)} = 
  \begin{bmatrix}
    K & O_{d,c} & -T_{\rr{n}} & -T_{\rr{t}} \\
  \end{bmatrix}
  ,   \quad
  \bi{w}_{2}^{(2)} = -\bi{p} . 
\end{align*}
Then we see that \eqref{P.linear.SOC.2.1} is reduced to 
\begin{align*}
  M_{22}^{(2)}  \bi{v} 
  + \bi{w}_{2}^{(2)} = \bi{0} . 
\end{align*}
Consequently, problem \eqref{P.linear.SOC.2} can be transformed into the 
the form in \eqref{P.linear.SOC.1} with 
\begin{alignat*}{3}
  & M_{11}= O_{4c,4c} , 
  &{\quad}& 
  &{\quad}& 
  \bi{w}_{1} = O_{4c,1} , \\
  & M_{21} = 
  \begin{bmatrix}
    I_{4c} \\ O_{d,4c} \\
  \end{bmatrix}
  , 
  &{\quad}& 
  M_{22} = 
  \begin{bmatrix}
    M_{22}^{(1)} \\
    M_{22}^{(2)} \\
  \end{bmatrix}
  , 
  &{\quad}& 
  \bi{w}_{2} = 
  \begin{bmatrix}
    \bi{w}_{2}^{(1)} \\
    \bi{w}_{2}^{(2)} \\
  \end{bmatrix}
  ,
\end{alignat*}
and $M_{12}$ in \eqref{eq:matrix.M.1.2}.

\section{Optimality condition of problem \eqref{P.approximation.3}}
\label{sec:optimality}

We can confirm equivalence of \eqref{eq:modified.complementarity} and 
\eqref{P.approximation.3} as follows. 

Observe that problem \eqref{P.approximation.3} is 
written in an explicit manner as follows: 
\begin{subequations}\label{P.approximation.appendix.1}%
  \begin{alignat}{3}
    & \MIN_{\bi{u} \in \Re^{d}} 
    &{\quad}& 
    \frac{1}{2} \Delta\bi{u}^{\top} \, K \, \Delta\bi{u} 
    - \bi{p}^{\top} \, \Delta\bi{u} 
    \\
    & \ST && 
    \begin{bmatrix}
      -\tilde{g}_{j} + \bi{t}_{\rr{n}j}^{\top} \, \Delta\bi{u} \\
      T_{\rr{t}j}^{\top} \, \Delta\bi{u} \\
    \end{bmatrix}
    \in F^{*} , 
    \quad j=1,\dots,c . 
  \end{alignat}
\end{subequations}
Since $3c < d$, problem \eqref{P.approximation.appendix.1} has 
an interior feasible solution. 
The Lagrangian of problem \eqref{P.approximation.appendix.1} is defined by 
\begin{align}
  L(\Delta\bi{u};\bi{r})  
  = 
  \begin{dcases*}
    \frac{1}{2} \Delta\bi{u}^{\top} \, K \, \Delta\bi{u} 
    - \bi{p}^{\top} \, \Delta\bi{u} \notag\\
    \qquad
    {}- \sum_{j=1}^{c} \left\langle  
    \begin{bmatrix}
      r_{\rr{n}j} \\ \bi{r}_{\rr{t}j} \\
    \end{bmatrix}
    , 
    \begin{bmatrix}
      -\tilde{g}_{j} + \bi{t}_{\rr{n}j}^{\top} \, \Delta\bi{u} \\
      T_{\rr{t}j}^{\top} \, \Delta\bi{u} \\
    \end{bmatrix}
    \right\rangle
    & if 
    $\begin{bmatrix}
       r_{\rr{n}j} \\ \bi{r}_{\rr{t}j} \\
     \end{bmatrix}
    \in F$ $(j=1,\dots,c)$, \\
    -\infty 
    & otherwise. 
  \end{dcases*}
\end{align}
Indeed, from the duality between $F$ and $F^{*}$ we see that problem 
\eqref{P.approximation.appendix.1} is equivalent to the following one: 
\begin{align*}
  \MIN_{\Delta\bi{u} \in \Re^{d}} 
  \quad \sup_{\bi{r} \in \Re^{3c}} L(\Delta\bi{u};\bi{r}) . 
\end{align*}
Here, 
$\displaystyle \sup_{\bi{r} \in \Re^{3c}} L(\Delta\bi{u};\bi{r})$ 
attains at a finite value if and only if 
\eqref{eq:modified.complementarity.2} is satisfied. 
Also, the stationarity condition of 
$L(\Delta\bi{u};\bi{r})$ with respect to $\Delta\bi{u}$ is 
\eqref{eq:modified.complementarity.1}. 
Thus, \eqref{eq:modified.complementarity} is a necessary and sufficient 
condition for optimality of problem \eqref{P.approximation.3}.

\end{document}